\documentclass[10pt]{amsart}

\usepackage{amsmath,amsfonts,amssymb}

\newcommand{\R}{\mathbb{R}}
\newcommand{\N}{\mathbb{N}}
\newcommand{\car}{{\raise0pt\hbox{{\LARGE $\chi$}}}}
\newcommand{\sg}{{\rm \; sign}}
\newcommand{\sgb}{{\rm \; {\bf sign}}}
\newcommand{\sop}{{\rm supp}\,}
\newcommand{\Div}{\hbox{\rm div\,}}
\newcommand{\e }{\epsilon }

\newcommand{\dis }{{\mathcal D}' }
\newcommand{\z }{{\bf z}}

\newenvironment{pf}{\noindent{\sc Proof}.\enspace}{\rule{2mm}{2mm}\medskip}

\newtheorem{Theorem}{Theorem}[section]

\newtheorem{Definition}[Theorem]{Definition}

\newtheorem{Proposition}[Theorem]{Proposition}

\newtheorem{Remark}[Theorem]{Remark}

\newcommand{\norma}[2]{\|#1\|_{\lower 4pt \hbox{$\scriptstyle #2$}}}

\newcommand{\tod}{{\rightharpoonup \,}}
\newcommand{\tode}{\overset{*}{\rightharpoonup\,}}
\newcommand{\h}{{\mathcal H}^{N-1}}
\renewcommand{\b}{{\mathcal B}}

\begin{document}

\title[Elliptic 1-Laplacian equations with dynamical boundary conditions]{Elliptic 1-Laplacian equations with dynamical boundary conditions}

\author[M. Latorre and  S. Segura de Le\'on]
{Marta Latorre and  Sergio Segura de Le\'on}

\address{M. Latorre: Departament d'An\`{a}lisi Matem\`atica,
Universitat de Val\`encia,
Dr. Moliner 50, 46100 Burjassot, Spain.
{\it E-mail address:} {\tt marta.latorre@uv.es }}

\address{S. Segura de Le\'on: Departament d'An\`{a}lisi Matem\`atica,
Universitat de Val\`encia,
Dr. Moliner 50, 46100 Burjassot, Spain.
{\it E-mail address:}  {\tt sergio.segura@uv.es }}

\thanks{}


\bigskip
\begin{abstract}
This paper is concerned with an evolution problem having an elliptic equation involving the $1$-Laplacian operator and a dynamical boundary condition. We apply nonlinear semigroup theory to obtain existence and uniqueness results as well as a comparison principle. Our main theorem shows that the solution we found is actually a strong solution. We also compare solutions with different data.
\end{abstract}

\maketitle


\section{Introduction}

In this paper we deal with existence and uniqueness for an evolution problem. It consists in an elliptic equation involving the $1$-Laplacian operator and a dynamical boundary condition, namely,
\begin{equation} \label{Dynamic-globalproblem}
\left\{
\begin{array}{ll}
\displaystyle \lambda u-\Div \left(\frac{Du}{|Du|}\right) =0 & \mbox{ in } (0,+\infty)\times \Omega \,, \\[5mm]
\displaystyle \omega_t+\left[\frac{Du}{|Du|},\nu\right] = g(t,x) & \mbox{ on } (0,+\infty) \times \partial \Omega \,, \\[5mm]
u=\omega & \mbox{ on } (0,+\infty) \times \partial \Omega \,, \\[3mm]
\omega(0, x)=\omega_0(x) & \mbox{ on } \partial \Omega \,;
\end{array}\right.
\end{equation}
where $\Omega$ is a bounded open set in $\R^N$ with smooth boundary $\partial\Omega$, $\lambda>0$, $\nu$ stands for the unit outward normal vector on $\partial\Omega$, $g\in L_{loc}^1(0,+\infty;L^2(\partial\Omega))$ and $\omega_0\in L^2(\partial \Omega)$.  Here, we have denoted by $\omega_t$ the distributional derivative of $\omega$ with respect to $t$. As far as we know, this is the first time that dynamical boundary conditions for the $1$-Laplacian are considered.

We point out that dynamical boundary conditions naturally occur in applications where there is a reaction term in the problem that concentrates in a small strip around the boundary of the domain, while in the interior there is no reaction and only diffusion matters. So, it appears in many mathematical models including heat transfer in a solid in contact with a moving fluid, in thermoelasticity, in biology, etc. This fact has given rise to many papers  (see \cite{AF, AIMT1, AIMT2, BDV, CR, E1, E2, FH, G, H, L, RS}) dealing with problems having dynamical boundary conditions, and mainly of those problems involving linear operators. The study of problems where an elliptic or parabolic equation occurs with this kind of boundary conditions is nowadays an active branch of research and we refer to \cite{FIK, GW, Gu, MMR, ST}  and references therein for recent papers.

The study of an evolution problem having an elliptic equation driven by the $p$-Laplacian (with $p>1$) and a dynamical boundary condition is due to \cite{AIMT1} (see also \cite{AIMT2}). To handle with that nonlinear problem, the authors define a completely accretive operator, apply the nonlinear semigroup theory to get a mild solution and finally, prove that this mild solution is actually a weak solution. Once their result is available, we may study problem \eqref{Dynamic-globalproblem} taking the solution corresponding to $p>1$ and letting $p$ go to $1$. Nevertheless, we are not able to pass to the limit and this approach remains an open problem. Furthermore, once a solution to our problem is obtained, we cannot prove that it is the limit of mild solutions to problems involving the $p$-Laplacian. What we need to prove the convergence would be a Modica type result on lower semicontinuity (see \cite[Proposition 1.2]{Mo}) for functionals depending on time.

Instead trying this approach, we adapt the method used in \cite{AIMT1} and apply the nonlinear semigroup theory (we refer to \cite{BCP} for a good introduction to this theory). Obviously, the singular features of the $1$-Laplacian do not allow us to follow every step. Among the special features verified by the $1$-Laplacian, we highlight that boundary conditions need not be satisfied in the sense of traces (we refer to \cite{ABCM2} for the Dirichlet problem, to \cite{MRST} for the Neumann problem as well as \cite{ABCM1} for the homogeneous Neumann for a related equation, and to \cite{MRS2} for the Robin  problem). This fact leads us to modify the procedure from the very beginning since it implies a change in the definition of the associated accretive operator.
Indeed, the translation of the operator studied in \cite{AIMT1} to our setting would be an operator $\mathfrak{B} \subset L^2(\partial\Omega)\times L^2(\partial\Omega)$  defined as follows:
\begin{Definition}
Let  $v, \omega\in L^2(\partial\Omega)$. Then
$v\in \mathfrak{B} (\omega)$ if there exists $u\in BV(\Omega)\cap L^2(\Omega)\cap L^2(\partial\Omega)$ such that $u\big|_{\partial\Omega}=\omega$ and it is a solution to the Neumann problem
\begin{equation*}
  \left\{\begin{array}{ll}
 \displaystyle \lambda u-\Div\left(\frac{Du}{|Du|}\right)=0\,,&\hbox{in }\Omega\,;\\[5mm]
  \displaystyle \left[\frac{Du}{|Du|}\cdot\nu\right]=v\,,&\hbox{on }\partial\Omega\,.
  \end{array}\right.
\end{equation*}
\end{Definition}
This is indeed a completely accretive operator but, unfortunately, we are not able to prove that it satisfies the range condition; thus the nonlinear semigroup theory cannot be applied. We turn out to define our operator for $v, \omega\in L^2(\partial\Omega)$ as $v\in \mathcal{B} (\omega)$ if  $v\in L^\infty(\partial\Omega)$, with $\|v\|_{L^\infty(\partial\Omega)}\le 1$, and there exists $u\in BV(\Omega)\cap L^2(\Omega)$ which is a solution to the Dirichlet problem with datum $\omega$ and it is also a solution of the Neumann problem with datum $v$ (see Definition \ref{def.operator-B} below).
Now, we do not know if this operator is completely accretive, we only prove that it is accretive in $L^2(\partial\Omega)$. Hence, we have not to expect that our solution holds every feature satisfied by solutions to problems driven by the $p$-Laplacian (for instance, we just choose initial data belonging to $L^2(\partial\Omega)$). Moreover, even when our solution satisfies the same property, the proof of this fact can be different, as can be checked in the comparison principle. Despite these difficulties, we obtain global existence and uniqueness of solution for every datum $\omega_0\in L^2(\partial\Omega)$ as well as a comparison principle. Furthermore, we prove that the solution we found is a strong solution in the sense that the problem holds for almost all $t>0$. We also analyze some related properties as the continuous dependence on data. Our main result is the following.

\begin{Theorem}
Let $\lambda >0$, and let $g \in L_{loc}^1(0,+\infty;L^2(\partial\Omega))$ and $\omega_0 \in L^2(\partial \Omega)$. There exists a unique global solution  $(u, \omega) $  to problem \eqref{Dynamic-globalproblem} in the sense of Definition \ref{def.sol.Dynamic-problem}. This solution satisfies
$u \in L_{loc}^2(0,+\infty;L^2(\Omega))\cap L_{loc}^\infty(0,+\infty;BV(\Omega))$ and $\omega \in C([0,+\infty[; L^2(\partial \Omega))\cap W_{loc}^{1,1}(0,+\infty;L^2(\partial \Omega))$.

Furthermore, the following estimates hold:
\begin{equation*}
  \|\omega\|_{L^\infty (0,T; L^2(\partial\Omega))}\le \|\omega_0\|_{L^2(\partial\Omega)}+\|g\|_{L^1 (0,T; L^2(\partial\Omega))}\quad\hbox{ for every }T>0\,,
\end{equation*}
\begin{equation*}
  \lambda\|u(t)\|_{L^2(\Omega)}^2+\|u(t)\|_{BV(\Omega)} \le  \|\omega(t)\|_{L^1(\partial\Omega)}\quad\hbox{ for almost all }t>0\,.
\end{equation*}
\end{Theorem}

The paper is organized into 5 sections. In Section 2, we introduce our notation and state the main features of functions of bounded variation, of $L^\infty$-divergence-measure vector fields  and the theory of nonlinear semigroups. Section 3 is devoted to obtain the mild solution to the associated abstract Cauchy problem, while in Section 4 we check that this mild solution is actually a strong solution to problem \eqref{Dynamic-globalproblem}. Finally, Section 5 deals with continuous dependence of data.

\section{Preliminaries}

In this section, we will present some useful results and the notation used in what follows.

Throughout this paper, $\Omega$ is an open bounded set in $\R^N$ which boundary $\partial \Omega$ is smooth. So there exists the outward normal unit vector $\nu(x)$ for $\h$-almost every $x\in\partial\Omega$, where $\h$ denotes the $(N-1)$-dimensional Hausdorff measure.

For every $k>0$, we define the truncation function as
\begin{equation*}
T_k(s) = \min\{ |s|,k \}\sg (s) \,, \quad s \in \R \,.
\end{equation*}

We will work with the usual Lebesgue and Sobolev spaces, denoted by $L^q(\Omega)$  and $W_0^{1,p}(\Omega)$, respectively (see for instance \cite{Br} or \cite{Ev}).
If $T>0$, the spaces $L^r(0,T;L^{q}(\Omega))$ are defined as follows:
$$
u\in L^r(0,T;L^{q}(\Omega))
$$
if $u:(0,T)\times\Omega\to\R$ is Lebesgue measurable and the integral
\begin{equation*}
\int_0^T\bigg(\int_\Omega
|u(t, x)|^qdx\bigg)^{\frac{r}{q}}dt
\end{equation*}
is finite.
It is clear that  for $q,r\ge 1$, the space $L^r(0,T;L^{q}(\Omega))$ is a Banach space equipped with the norm
$$
||u||_{L^r(0,T;L^{q}(\Omega))}=\Big(\int_0^T\Big(\int_\Omega |u(t, x)|^qdx\Big)^{\frac{r}{q}}dt\Big)^{\frac 1r}.
$$
In a similar way we define the space $L^r(0,T;W^{1,q}_0(\Omega))$ or $W^{1,r}(0,T;L^q(\Omega))$. We refer to \cite{Ev} for more details.

Given a Banach function space $X$, recall that $u\in  L^r(0,T;X)$ implies that \mbox{$u(t)\in X$} for almost all $t\in \, ]0,T[$.
Moreover, instead of writing ``$u\in L^r(0,T;X)$ for every $T>0$'', we shall write $u\in L^r_{\rm loc}(0,+\infty;X)$.
Moreover, if $\mathcal I$ is a real interval, then $C(\mathcal I; X)$ stands for the space of all continuous functions from $\mathcal I$ into $X$.

\subsection{Functions of bounded variation}

The natural energy space to study our problem is the space of functions of bounded variation, denoted by $BV(\Omega)$. We say that a function \mbox{$u\,:\,\Omega\to\R$} belongs to $BV(\Omega)$ if $u \in L^1(\Omega)$ and its gradient in the sense of distributions $Du$ is a Radon measure with finite total variation.
The norm associated to this space is given by
\begin{equation*}
\|u\| =\int_\Omega |u| \,dx +\int_\Omega |Du| \,.
\end{equation*}
We recall that every function of bounded variation has a trace on the boundary, so that we may write $u\big|_{\partial\Omega}$.
Moreover, there exists a bounded linear operator $BV(\Omega)\hookrightarrow L^1(\partial\Omega)$ which is also onto. As a consequence, an equivalent norm on $BV(\Omega)$ can be defined:
\begin{equation*}
\|u\|_{BV(\Omega)} =\int_{\partial\Omega} |u| \,d\h +\int_\Omega |Du| \,.
\end{equation*}
We will often use this norm in what follows.

Throughout this paper we have to use the lower semicontinuity of some functionals  defined on $BV(\Omega)$ with respect to the convergence in $L^1(\Omega)$. The result we will apply is stated as follows.

\begin{Theorem}\label{teo.semicont-inferior}
Let $\omega \in L^1(\partial \Omega)$ and let $\varphi \in C_0^1(\Omega)$ with $\varphi \ge0$. If the sequence $\{u_n\} \subseteq BV(\Omega)$ converges to $u$ in $L^1(\Omega)$, then the following inequalities hold
\begin{equation*}
\int_\Omega |Du| + \int_{\partial\Omega} |u-\omega| \,d\h  \le \liminf_{n \to \infty} \int_\Omega |Du_n| + \int_{\partial\Omega} |u_n-\omega| \,d\h \,
\end{equation*}
and
\begin{equation*}
\int_\Omega \varphi \, |Du| \le \liminf_{n \to \infty} \int_\Omega \varphi \, |Du_n| \,.
\end{equation*}
\end{Theorem}

%

For further information about functions of bounded variation we refer to \cite{AFP}, \cite{EG} and \cite{Zi}.

\subsection{Green's formula}
Following \cite{ABCM1}, the quotient $\dfrac{Du}{|Du|}$ in our equation makes sense through a vector field $\z\in L^\infty(\Omega;\R^N)$ satisfying two conditions: (i) $\|\z\|_\infty\le1$ and (ii) the dot product of $\z$ and $Du$ is equal to $|Du|$. The validity of this dot product lies on the Anzellotti theory (see \cite{An}). Consider $\z \in L^\infty (\Omega;\R^N)$ such that $\Div \z \in L^2(\Omega)$ and $u \in BV(\Omega)\cap L^2(\Omega)$ and define the functional
\begin{equation*}
\langle (\z, Du),\varphi \rangle = - \int_\Omega u \,\varphi \, \Div \z - \int_\Omega u \, \z \cdot \nabla \varphi \,dx \,,
\end{equation*}
for every $\varphi \in C_0^\infty (\Omega)$. This distribution turns out to be a Radon measure such that its total variation satisfies
\begin{equation*}
  |(\z, Du)|\le \|\z\|_\infty|Du|\quad\hbox{as measures.}
\end{equation*}


Due to the Anzellotti theory, a definition of a weak trace on $\partial \Omega$ of the normal component of $\z$ is given, it is denoted by $[\z,\nu]$ and it satisfies the inequality $\|[\z,\nu]\|_{L^\infty(\partial\Omega)}\le \|\z\|_{L^\infty(\Omega)}$. Moreover, a Green's formula involving all these elements holds:
\begin{Theorem}
If $\z \in L^\infty (\Omega;\R^N)$ satisfies $\Div \z \in L^2(\Omega)$ and $u\in BV(\Omega)\cap L^2(\Omega)$, then it holds
\begin{equation*}
\int_{\Omega} u \,\Div\z + \int_{\Omega} (\z, Du) =
\int_{\partial \Omega} u\,[\z, \nu] \ d\mathcal H^{N-1}\,.
\end{equation*}
\end{Theorem}

Although we usually take the above assumptions, we point out that $(\z, Du)$ can be defined for other pairings; for instance, $\Div \z \in L^N(\Omega)$ and $u \in BV(\Omega)$ or $\Div \z \in L^1(\Omega)$ and $u \in BV(\Omega)\cap L^\infty(\Omega)$. In every case, the results stated above also hold true.

\subsection{Mild solutions}
In this subsection we will present some definitions and results concerning mild solutions.

Let $X$ be a Banach space and let $\mathcal{P}(X)$ be the collection of all subsets of $X$. Every mapping $\mathcal{A} : X \to \mathcal{P}(X)$ will be called an operator in $X$.
\begin{Definition}
An operator $\mathcal{A}:X \to \mathcal{P}(X)$ is said to be accretive if
\begin{equation*}
\| v-\widehat{v} + \alpha(\omega - \widehat{\omega})\|_X \ge \|v-\widehat{v}\|_X \,,
\end{equation*}
whenever $\alpha \ge 0$, and $v \in \mathcal{A}(\omega)$ and $\widehat{v} \in \mathcal{A}(\widehat{\omega})$. When $X$ is a Hilbert space, the operator $\mathcal{A}$ is accretive if and only if it is monotone, that is,
\begin{equation*}
\langle v-\widehat{v},\omega-\widehat{\omega}\rangle \ge 0 \,,
\end{equation*}
for every $v \in \mathcal{A}(\omega)$ and $\widehat{v} \in \mathcal{A}(\widehat{\omega})$.
\end{Definition}

\begin{Definition}
An operator $\mathcal{A}:X \to \mathcal{P}(X)$ is $m$-accretive if it is accretive and $R(I+\e \mathcal{A})=X$ for all  $\e >0$.
\end{Definition}

We next introduce the notion of mild solution to the abstract Cauchy problem
\begin{equation}\label{ec.mild-sol}
\left\{\begin{array}{ll}
\omega_t + \mathcal{A}(\omega) \ni g\,,\\[3mm]
\omega(0)=\omega_0\,,
\end{array}\right.
\end{equation}
where $g\in L^1_{loc}(0,+\infty; X)$ and $\omega_0\in X$.

\begin{Definition}\label{def.discretization}
Fix $T>0$.
If $t_0<t_1 < \dots <t_n$ satisfy
\begin{equation*}\begin{array}{ll}
0\le t_0<\e \,, &\\[2mm]
t_i-t_{i-1} < \e \quad &\mbox{ for }\quad i=1,2, \dots, n \,,\\[2mm]
0\le T-t_n<\e\,, &
\end{array}
\end{equation*}
and $g_1, g_2, \dots,g_n$ is a finite sequence in $X$ such that
\begin{equation*}
\sum_{i=1}^n \int_{t_{i-1}}^{t_i} \| g(s)-g_i\|_X \,ds< \e \,,
\end{equation*}
then the system
\begin{equation}\label{ec.discretizada}
\frac{\omega_{i}-\omega_{i-1}}{t_i-t_{i-1}} + \mathcal{A}(\omega_i) \ni g_i \qquad\mbox{ for }\quad i=1,2,\dots,n \,,
\end{equation}
is called an $\e$-discretization of \eqref{ec.mild-sol} on $[0,T]$.

We say that a function $\omega_\e:[t_0,t_n] \to X$ is a solution to this $\e$-discretization  if $\omega_\e$ is a piecewise constant function such that \mbox{$\omega_\e(t_0)=\omega_0$}, $\omega_\e(t)=\omega_i$ on $]t_{i-1},t_i]$ for $i=1,2,\dots,n$, and system \eqref{ec.discretizada} holds.
\end{Definition}

\begin{Remark}
\rm
Definition \ref{def.discretization} is based on the possibility of approximating any function $g\in L^1(0,T; X)$ by steps functions $\sum_{i=1}^n g_i \chi_{]t_{i-1}, t_i]}$. We point out that this approximation can be taken in such way that $g_i=g(t_i)$, being $t_i$ a Lebesgue point of $g$ for $i=1,\dots ,n$ (see \cite[Proposition 1.5]{BCP}).
\end{Remark}

\begin{Definition}
Fix $T>0$ and let $g \in L^1(0,T;X)$. A mild solution of the abstract Cauchy problem \eqref{ec.mild-sol} on $[0,T]$ is a function $\omega \in C([0,T];X)$ such that, for every $\e >0$, there exists an $\e$-discretization of \eqref{ec.mild-sol} on $[0,T]$ which has a solution $\omega_\e$ satisfying
\begin{equation*}
\|\omega(t)-\omega_\e(t)\|_X < \e \qquad \mbox{ for all } \quad t \in [0,T] \,.
\end{equation*}
\end{Definition}

\begin{Definition}
Let $g \in L^1_{loc}(0,+\infty;X)$. A mild solution of problem \eqref{ec.mild-sol} on $[0,+\infty[$ is a function $\omega \in C([0,+\infty[; X)$ whose restriction to each  subinterval $[0,T]$ of $[0,+\infty[$ is a mild solution on $[0,T]$.
\end{Definition}

\begin{Remark}\label{convergencia.mild-solutions}\rm
From the definition of mild solution one deduces that solutions to discretizations satisfy
\begin{equation*}
\omega_\e \to \omega \qquad \mbox{ in } L^\infty([0,T];X) \,,
\end{equation*}
for every $T>0$.
\end{Remark}

\begin{Theorem}\label{teo.existencia.mild.solutions}
Let $\mathcal{A}$ be an $m$-accretive operator in $X$. Consider $\omega_0 \in \overline{D(\mathcal{A})}$ and $g \in L^1_{loc}([0,+\infty[;X)$. Then problem \eqref{ec.mild-sol}  has a unique mild solution $\omega$ on $[0,+\infty[$.
\end{Theorem}

A final definition is in order.

\begin{Definition}
Fix  $T>0$ and let $g \in L^1(0,T;X)$. A strong solution of problem \eqref{ec.mild-sol} on $[0,T]$ is an absolutely continuous function $\omega \>:\>[0,T]\to X$ which is differentiable almost everywhere on $[0,T]$ and satisfies $\omega_t(t) + \mathcal{A}(\omega(t)) \ni g(t)$ for almost all $t\in [0,T]$.
\end{Definition}

We point out that every strong solution is a mild solution (see \cite[Theorem 1.4]{BCP}), but the converse does not hold.

For further information about mild solutions and semigroups on Banach spaces we refer to \cite{BCP} (and to \cite{Br0} for semigroups on Hilbert spaces).

\section{Existence of mild solutions}
Let $T>0$ and consider the problem
\begin{equation} \label{Dynamic-problem}
\left\{
\begin{array}{ll}
\displaystyle \lambda u-\Div \left(\frac{Du}{|Du|}\right) =0 & \mbox{ in } (0,T)\times \Omega \,, \\[5mm]
\displaystyle \omega_t+\left[\frac{Du}{|Du|},\nu\right] = g(t,x) & \mbox{ on } (0,T) \times \partial \Omega \,, \\[5mm]
u=\omega & \mbox{ on } (0,T) \times \partial \Omega \,, \\[3mm]
\omega(0, x)=\omega_0(x) & \mbox{ on } \partial \Omega \,;
\end{array}\right.
\end{equation}

As we have already mentioned, we want to define an accretive operator in $L^2(\partial\Omega)$ to apply the semigroup theory and then get a mild solution. Afterwards, using this mild solution we will obtain a strong solution to problem \eqref{Dynamic-problem}.
\begin{Remark}\rm
We point out that our operator will be defined on the boundary, and so our mild solution is $\omega$, while $u$ appearing in problem \eqref{Dynamic-problem} is just the corresponding auxiliary function. Nevertheless, this auxiliary function $u$ is univocally determined by $\omega$, since solutions to the Dirichlet problem for equation $\lambda u-\Div \left(\frac{Du}{|Du|}\right) =0$ are unique (see \cite{ABCM2}).
\end{Remark}

We start defining the operator $\b$ in the space $L^2(\partial \Omega)$.
\begin{Definition}\label{def.operator-B}
We say that $v \in \b(\omega)$ if $\omega$ belongs to $L^2(\partial \Omega)$ and $v$ belongs to $L^\infty(\partial\Omega)$, with $\|v\|_{L^\infty(\partial\Omega)}\le1$, and there exist a function $u \in BV(\Omega)\cap L^2(\Omega)$ and a vector field $\z \in L^\infty(\Omega;\R^N)$ with $\|\z\|_{L^\infty(\Omega)} \le 1$ such that
\begin{equation*}
\begin{array}{llccc}
(i) & \lambda u - \Div \z = 0 \qquad \mbox{ in } \dis(\Omega) \,, &\qquad&\qquad&\qquad\\[3mm]
(ii) & (\z,Du)=|Du| \qquad \mbox{ as measures in } \Omega\,, &\qquad&\qquad&\qquad\\[3mm]
(iii) &  [\z,\nu]=v \qquad \mathcal H^{N-1}\mbox{-a.e. on } \partial\Omega \,, &\qquad&\qquad&\qquad\\[3mm]
(iv) &  [\z,\nu]\in \sg(\omega-u) \qquad \mathcal H^{N-1}\mbox{-a.e. on } \partial \Omega \,. &\qquad&\qquad&\qquad\
\end{array}
\end{equation*}
\end{Definition}
Moreover, using Green's theorem and since conditions $(i)$ and $(iii)$ hold, we may deduce the following variational formulation:
\begin{equation*}
\lambda\int_\Omega u \, \varphi \,dx + \int_\Omega (\z,D\varphi) = \int_{\partial \Omega} v \, \varphi \, d\h \,,
\end{equation*}
for every test function $\varphi \in BV(\Omega)\cap L^2(\Omega)$. Notice that function $v \in L^\infty(\partial\Omega)$ and $\varphi\big|_{\partial\Omega}\in L^1(\partial\Omega)$, so that the last integral is well-defined.

In other words, we say that $v \in \b(\omega)$ if there exists $u \in BV(\Omega)\cap L^2(\Omega)$ such that $u$ is a solution to equation
\begin{equation}\label{ecuacion}
\lambda u-\Div\left(\frac{Du}{|Du|}\right)=0  \quad \mbox{ in } \Omega \,,
\end{equation}
with the Dirichlet boundary condition:
\begin{equation}\label{bound.cond.Dirichlet}
u = \omega \quad\mbox{ on } \partial \Omega \,,
\end{equation}
and it is also a solution to equation \eqref{ecuacion} with Neumann boundary condition
\begin{equation}\label{bound.cond.Neumann}
\left[\dfrac{Du}{|Du|},\nu\right] = v \quad \mbox{ on } \partial \Omega \,.
\end{equation}

From another point of view, operator $\b$ can be written as $v \in \b (\omega)$ if $v \in L^\infty(\partial\Omega)$ satisfies
\begin{itemize}
\item[$(i)$] $\|v\|_{L^\infty(\partial\Omega)} \le 1$,
\item[$(ii)$] $v \in \sg(\omega-u)$ where $u$ is the solution to \eqref{ecuacion} with boundary condition \eqref{bound.cond.Neumann}.
\end{itemize}
\subsection{Associated Robin problem}
Now, we analyze the Robin problem for \eqref{ecuacion}, to this end we follow \cite{MRS2}. For $\beta>0$, we consider the boundary condition:
\begin{equation}\label{bound.cond.Robin}
\beta u + \left[\frac{Du}{|Du|}, \nu\right] = g \quad \mbox{ on } \partial \Omega \,.
\end{equation}

\begin{Definition}
Let $g \in L^2(\partial\Omega)$, we say that $u \in BV(\Omega) \cap L^2(\Omega)$ is a weak solution to Robin problem \eqref{bound.cond.Robin} for equation \eqref{ecuacion} if there exists a vector field $\z \in L^\infty(\Omega;\R^N)$ with $\|\z\|_\infty \le 1$ such that
\begin{equation*}
\begin{array}{ll}
(i) & \lambda u-\Div\z=0 \qquad \mbox{ in } \mathcal D'(\Omega)\,,\\[3mm]
(ii) & (\z,Du)=|Du| \qquad \mbox{ as measures in } \Omega\,, \\[3mm]
(iii) & T_1(\beta u -g) = - [\z,\nu] \qquad \mathcal H^{N-1}\mbox{-a.e. on } \partial \Omega \,. \\[3mm]
\end{array}
\end{equation*}
\end{Definition}
As a consequence of Green's formula, the following variational formulation holds:
\begin{equation}\label{var1}
  \displaystyle \lambda\int_\Omega u \, \varphi \,dx + \int_\Omega (\z,D\varphi ) + \int_{\partial \Omega} T_1(\beta u -g) \, \varphi \,d\h = 0 \,,
\end{equation}
for every $\varphi \in BV(\Omega)\cap L^2(\Omega)$.

\begin{Remark}\label{remark.sol-Dirichlet-problem}\rm
Every solution to equation \eqref{ecuacion} with the Robin boundary condition \eqref{bound.cond.Robin} is also a solution to the same equation but with Dirichlet boundary condition \eqref{bound.cond.Dirichlet} for $\omega$ satisfying $T_1(\beta u-g)=\beta \omega -g$ (see \cite[Proposition 2.13]{MRS2}). Using this function, \eqref{var1} becomes
\begin{equation*}
  \displaystyle \lambda\int_\Omega u \, \varphi \,dx + \int_\Omega (\z,D\varphi ) + \int_{\partial \Omega} (\beta \omega -g) \, \varphi \,d\h = 0 \,,
\end{equation*}
for every $\varphi \in BV(\Omega)\cap L^2(\Omega)$.
\end{Remark}

\begin{Remark}\label{comp2}\rm
Consider $g_1, g_2\in L^2(\partial\Omega)$ and let $u_i\in BV(\Omega)\cap L^2(\Omega)$ be the corresponding solutions to the Robin problem. Denote by $\z_i\in L^\infty(\Omega; \R^N)$  the associated vector fields and by $\omega_i$ the functions satisfying $T_1(\beta u_i-g_i)=\beta \omega_i -g_i$, for $i=1,2$.  Now, we can prove that
$g_1\le g_2$ on $\partial\Omega$ implies $u_1\le u_2$ in $\Omega$ and $\omega_1\le \omega_2$ on $\partial\Omega$.  It is enough to  take $\varphi=(u_1-u_2)^+$ as test function in the respective variational formulations and perform straightforward manipulations to obtain
\begin{equation}\label{comp1}
  \lambda\int_\Omega\big[(u_1-u_2)^+\big]^2\,dx\le \int_{\partial\Omega}(T_1(g_1-\beta u_1)-T_1(g_2-\beta u_2))(u_1-u_2)^+\, d\h\,.
\end{equation}
Note that, on the set $\{u_1\big|_{\partial\Omega}\ge u_2\big|_{\partial\Omega}\}$, the assumption $g_1\le g_2$ implies
\[
T_1(g_1-\beta u_1)-T_1(g_2-\beta u_2)\le0\,.
\]
Thus, the right hand side of \eqref{comp1} is nonnpositive and so  $(u_1-u_2)^+$ vanishes in $\Omega$. Moreover,
\[
\beta\omega_1=g_1-T_1(g_1-\beta u_1)\le g_2-T_1(g_2-\beta u_1)\le g_2-T_1(g_2-\beta u_2)=\beta\omega_2\,,
\]
so that $\omega_1\le \omega_2$ on $\partial\Omega$.
\end{Remark}

\subsection{Main properties of $\b$}
In this subsection, we will see the main properties of operator $\b$ that lead to a mild solution of problem \eqref{Dynamic-problem}.
We begin by showing that our operator is accretive.

\begin{Theorem}
The operator $\b$ given in Definition \ref{def.operator-B} is accretive in $L^2(\partial\Omega)$.
\end{Theorem}
\begin{pf}
Since $L^2(\partial\Omega)$ is a Hilbert space, we just have to prove that $\b$ is monotone.
\\
Let $v_i \in \b(\omega_i)$ for $i=1,2$. We will show that
\begin{equation*}
\int_{\partial\Omega} (v_1-v_2)(\omega_1-\omega_2) \,d\h \ge 0 \,.
\end{equation*}

Given $v_i \in \b(\omega_i)$, we may find functions $u_i \in BV(\Omega)\cap L^2(\Omega)$ and vector fields $\z_i \in L^\infty(\Omega;\R^N)$ with $\|\z_i\|_{L^\infty(\Omega)} \le 1$ such that
\begin{equation*}
\begin{array}{llccc}
(i) & \lambda u_i - \Div \z_i = 0 \qquad \mbox{ in } \dis(\Omega) \,, &\qquad&\qquad&\qquad\\[3mm]
(ii) & (\z_i,Du_i)=|Du_i| \qquad \mbox{ as measures in } \Omega\,, &\qquad&\qquad&\qquad\\[3mm]
(iii) &  [\z_i,\nu_i]=v_i \qquad \mbox{ on } \partial\Omega \,, &\qquad&\qquad&\qquad\\[3mm]
(iv) &  [\z_i,\nu_i]\in \sg(\omega_i-u_i) \qquad \mbox{ on } \partial \Omega \,, &\qquad&\qquad&\qquad \\[3mm]
(v) & \displaystyle \lambda\int_\Omega u_i \, \varphi \,dx + \int_\Omega (\z_i,D\varphi)  = \int_{\partial \Omega} v_i \, \varphi \, d\h \,, &\qquad&\qquad&\qquad\\[3mm]
\end{array}
\end{equation*}
for every $\varphi \in BV(\Omega)\cap L^2(\Omega)$ and for $i=1,2$. Taking $u_1-u_2$ as a test function in $(v)$ for both $i=1,2$ and subtracting one from the other, we get
\begin{align*}
\lambda & \int_\Omega (u_1-u_2)^2 \,dx +\int_\Omega \left[|Du_1|-(\z_2,Du_1)+|Du_2|-(\z_1,Du_2)\right]
\\ & = \int_{\partial \Omega} (v_1-v_2) (u_1-u_2) \, d\h \,.
\end{align*}
Since the left hand side is positive (note that $(z_i,Du_j)\le |Du_j|$ for $i,j=1,2$), we deduce that
\begin{align}\label{acretivo-ec-1}
0 \le & \int_{\partial \Omega} (v_1-v_2) (u_1-u_2) \, d\h=\int_{\partial \Omega} (v_1-v_2) (\omega_1-\omega_2) \, d\h
\\ \nonumber & +\int_{\partial \Omega} (v_1-v_2) (u_1-\omega_1) \, d\h +\int_{\partial \Omega} (v_2-v_1) (u_2-\omega_2) \, d\h \,.
\end{align}
On the one hand, using conditions $(iii)$ and $(iv)$ and that $\|v_i\|_\infty \le 1$, it holds
\begin{align*}
&\int_{\partial \Omega} (v_1-v_2) (u_1-\omega_1) \, d\h
\\ &\; = \int_{\partial \Omega} v_1(u_1-\omega_1) \, d\h - \int_{\partial \Omega} v_2 (u_1-\omega_1) \, d\h
\\ &\; =-\int_{\partial \Omega} |u_1-\omega_1| \, d\h - \int_{\partial \Omega} v_2 (u_1-\omega_1) \, d\h
\\ &\; =-\int_{\partial \Omega} \left(|u_1-\omega_1| +  v_2 (u_1-\omega_1)\right) \, d\h \le 0 \,,
\end{align*}
and similarly
\begin{equation*}
\int_{\partial \Omega} (v_2-v_1) (u_2-\omega_2) \, d\h = - \int_{\partial \Omega}\left( |u_2-\omega_2|  + v_1(u_2-\omega_2) \right)\, d\h \le 0 \,.
\end{equation*}

Therefore, using \eqref{acretivo-ec-1} we conclude that
\begin{equation*}
0 \le  \int_{\partial \Omega} (v_1-v_2) (u_1-u_2) \, d\h  \le \int_{\partial \Omega} (v_1-v_2) (\omega_1-\omega_2) \, d\h \,.
\end{equation*}
\end{pf}

\begin{Proposition}\label{prop.rango}
The operator $\b$ given in Definition \ref{def.operator-B} is $m$-accretive in $L^2(\partial\Omega)$.
\end{Proposition}
\begin{pf}
Denoting by $I$ the  identity operator in $L^2(\partial\Omega)$, we just have to prove
\begin{equation*}
R(I+\e\b)=L^2(\partial\Omega) \qquad \mbox{ for every }\quad \e>0 \,.
\end{equation*}
Given $\e>0$, it is enough to see that $L^2(\partial \Omega) \subseteq R(I+\e \b)$.

For every $g \in L^2(\partial\Omega)$, we will show that there exists $\omega \in L^2(\partial \Omega)$ such that $g \in \omega + \e \b(\omega)$. That is, we will see that $\dfrac{1}{\e}g - \dfrac{1}{\e} \omega \in \b(\omega)$.

We consider the following Robin problem
\begin{equation*}
\left\{
\begin{array}{ll}
\displaystyle \lambda u-\Div\left(\frac{Du}{|Du|}\right)=0 & \mbox{ in } \Omega \,, \\[5mm]
\displaystyle \dfrac{1}{\e}u+ \left[\frac{Du}{|Du|},\nu\right] =\dfrac{1}{\e}g & \mbox{ on } \partial \Omega \,.
\end{array}
\right.
\end{equation*}
Applying \cite[Theorem 1.1]{MRS2}, there exist a solution $u \in BV(\Omega)\cap L^2(\Omega)$, a vector field $\z \in L^\infty(\Omega;\R^N)$ with $\|\z\|_{L^\infty(\Omega)}\le 1$ and a function $\omega \in L^2(\partial \Omega)$ such that
\begin{equation*}
(\z,Du) =|Du| \quad \mbox{ as measures in } \Omega \,,
\end{equation*}
\begin{equation*}
[\z,\nu]=T_1\Big(\frac{1}{\e}g-\frac{1}{\e} u\Big)=\frac{1}{\e}g-\frac{1}{\e} \omega \,,
\end{equation*}
and so
\begin{equation*}
\frac{1}{\e}g-\frac{1}{\e} \omega\in L^\infty(\partial\Omega)\quad\hbox{with}\quad \Big\|\frac{1}{\e}g-\frac{1}{\e} \omega\Big\|_{L^\infty(\partial\Omega)}\le 1 \,.
\end{equation*}

In addition, $u$ is also a solution to the Dirichlet problem
\begin{equation*}
\left\{
\begin{array}{ll}
\lambda u-\Div\left(\frac{Du}{|Du|}\right)=0 & \mbox{ in } \Omega \,, \\
u = \omega & \mbox{ on } \partial \Omega \,,
\end{array}
\right.
\end{equation*}
(see Remark \ref{remark.sol-Dirichlet-problem}). Therefore, it also holds
\begin{equation*}
[\z,\nu] \in \sg (\omega-u) \,.
\end{equation*}
Thus, $\dfrac{1}{\e}g-\dfrac{1}{\e}\omega \in \b(\omega)$.
\end{pf}

\begin{Remark}\label{comp3}\rm
Proposition \ref{prop.rango} guarantees the existence of the resolvent
$$(I+\epsilon \b)^{-1}\>:\> L^2(\partial\Omega)\to L^2(\partial\Omega)$$
 for every $\epsilon>0$.
Taking into account Remark \ref{comp2}, we deduce that it is an order preserving operator.
\end{Remark}

\begin{Proposition}
Let $\b$ be the operator given in Definition \ref{def.operator-B}. Then, it holds
\begin{equation*}
L^2(\partial \Omega)=\overline{D(\b)} \,.
\end{equation*}
\end{Proposition}
\begin{pf}
We just have to prove that $L^2(\partial \Omega) \subseteq \overline{D(\b)}$. We begin by taking $g$ to be a function in $L^\infty(\partial\Omega)$. Given $n \in \N$, by Theorem \ref{prop.rango}, we know that $g \in R(I+\dfrac{1}{n}\b)$. Then, there exists $\omega_n \in L^2(\partial \Omega)$ such that $g \in \omega_n+\dfrac{1}{n}\b(\omega_n)$. That is, $n(g-\omega_n)\in \b(\omega_n)$. Therefore, there exist $u_n \in BV(\Omega)\cap L^2(\Omega)$ and a vector field $\z_n \in L^\infty(\Omega;\R^N)$ with $\|\z_n\|_{L^\infty(\Omega)}\le 1$ such that
\begin{align*}
 & (\z_n,Du_n)=|Du_n| \qquad \mbox{ as measures in } \Omega\,,
\\ & [\z_n,\nu]=n(g-\omega_n) \qquad \mbox{ on } \partial\Omega \,,
\\ & [\z_n,\nu]\in \sg(\omega_n-u_n) \qquad \mbox{ on } \partial \Omega \,,
\end{align*}
and
\begin{equation}\label{last}
\lambda\int_\Omega u_n \, \varphi \,dx + \int_\Omega (\z_n,D\varphi) = \int_{\partial \Omega} n(g-\omega_n) \, \varphi \, d\h \,,
\end{equation}
for every $\varphi \in BV(\Omega)\cap L^2(\Omega)$.

Since $g \in L^\infty(\partial\Omega)$, we have $g=v|_{\partial\Omega}$ for some $v \in W^{1,1}(\Omega)\cap L^\infty(\Omega)$ (see \cite{Ga}),
and we use $v-u_n$ as a test function in \eqref{last} to get
\begin{equation}\label{eq-1}
\lambda\int_\Omega u_n(v-u_n) \,dx + \int_\Omega(\z_n,D(v-u_n)) = \int_{\partial\Omega}n(g-\omega_n)(g-u_n) \,d\h \,.
\end{equation}
Observe that, since $n(g-\omega_n)\in \sg (\omega_n-u_n)$, we also have
\begin{equation}\label{eq-2}
\int_{\partial\Omega}n(g-\omega_n)(g-u_n)\,d\h  = n\int_{\partial\Omega}(g-\omega_n)^2\,d\h +\int_{\partial\Omega}|\omega_n-u_n|\,d\h \,.
\end{equation}
Joining now equations \eqref{eq-1} and \eqref{eq-2} we get
\begin{align*}
\lambda & \int_\Omega u_n v \,dx + \int_\Omega(\z_n,Dv) - \int_\Omega|Du_n|
\\ & = n\int_{\partial\Omega}(g-\omega_n)^2\,d\h +\int_{\partial\Omega}|\omega_n-u_n|\,d\h +\lambda\int_\Omega u_n^2 \,dx \,,
\end{align*}
and so it follows that
\begin{equation*}
\lambda\int_\Omega u_n v \,dx + \int_\Omega(\z_n,Dv)\ge n\int_{\partial\Omega}(g-\omega_n)^2\,d\h  +\lambda\int_\Omega u_n^2 \,dx \,.
\end{equation*}
Then, using Young's inequality and the fact that $(\z_n,Dv)\le |(\z_n,Dv)| \le |Dv|$ we obtain
\begin{equation*}
 n\int_{\partial\Omega}(g-\omega_n)^2\,d\h  +\lambda\int_\Omega u_n^2 \,dx  \le \frac{\lambda}{2}\int_\Omega u_n^2\,dx+\frac{\lambda}{2} \int_{\Omega} v^2 \,dx + \int_\Omega |Dv| \,.
\end{equation*}
Thus, simplifying,
\begin{equation*}
 n\int_{\partial\Omega}(g-\omega_n)^2\,d\h + \dfrac{\lambda}{2}\int_\Omega u_n^2 \,dx  \le \frac{\lambda}{2} \int_{\Omega} v^2 \,dx + \int_\Omega |Dv| \,,
\end{equation*}
and it yields
\begin{equation*}
\int_{\partial\Omega}(g-\omega_n)^2\,d\h \le \dfrac{1}{n} \left(\frac{\lambda}{2} \int_{\Omega} v^2 \,dx + \int_\Omega |Dv| \right)\,.
\end{equation*}
Finally, since the right-hand side goes to 0 as $n \to \infty$ we deduce that $\omega_n \to g$ in $L^2(\partial \Omega)$ and then, $g \in \overline{D(\b)}$.

Now, let $g \in L^2(\partial \Omega)$.
We already know that each truncation $T_k(g) \in\overline{D(\b)}$ and $T_k(g) \to g$ in $L^2(\partial \Omega)$ when $k$ goes to $+\infty$. Therefore, $g \in \overline{D(\b)}$.
\end{pf}

Using the previous results, the main theorem of this subsection can be obtained applying Theorem \ref{teo.existencia.mild.solutions}.
\begin{Theorem}\label{ex-mild}
Let $g \in L^1(0,T;L^2(\partial\Omega))$ and let $\omega_0\in L^2(\partial\Omega)$. Then there exists a unique mild solution to the abstract Cauchy problem $\omega_t+\b(\omega)\ni g$, $\omega(0)=\omega_0$ on $[0,T]$.
\end{Theorem}

\begin{Remark}\rm
 Some remarks concerning the limiting case $\lambda=0$ are in order. In this case, the definition of operator $\b$ must be modified, now the auxiliary function $u$ belongs to $BV(\Omega)$ (but, in general, not to $L^2(\Omega)$). Furthermore, now the definition of $(\z, Du)$ depends on the duality $\Div\z\in L^N(\Omega)$ and $u\in BV(\Omega)$. We point out that all the results proved in this section hold.

 Nevertheless, this auxiliary function $u$ is not longer determined by $\omega$ (see \cite{MRS1} for examples of nonuniqueness of the Dirichlet problem for the 1-Laplacian) and, moreover, the arguments of the next section does not work. Hence, we may prove that a mild solution exists, but we are not able to see that it is actually a strong solution.
\end{Remark}

\subsection{Comparison principle}

In this subsection, we will compare two mild solutions when their data are ordered.

\begin{Theorem}\label{compar}
Let $g^1, g^2 \in L^1(0,T;L^2(\partial\Omega))$ and let $\omega_0^1, \omega_0^2 \in L^2(\partial\Omega)$. Denote by $\omega^k \in C([0,T];L^2(\partial\Omega))$ the mild solution corresponding to data $g^k$ and  $\omega_0^k$, $k=1,2$.

If $g^1(t,x)\le g^2(t,x)$ for almost all $(t, x)\in (0,T)\times \partial\Omega$ and $\omega_0^1(x)\le \omega_0^2(x)$ for almost all $x\in \partial\Omega$, then the solutions to every $\e$-discretization satisfy $\omega_\e^1(t,x) \le \omega_\e^2(t,x)$ as well as the corresponding auxiliary functions $u_\e^1(t, x)\le u_\e^2(t, x)$.
As a consequence, $\omega^1(t,x)\le \omega^2(t,x)$ for almost all $(t, x)\in (0,T)\times \partial\Omega$.
\end{Theorem}

\begin{pf} Given $\epsilon>0$, consider an $\epsilon$-discretization of \eqref{ec.mild-sol} for data $g^k$ and  $\omega_0^k$. Observe that splitting the subintervals if necessary, we may take the same partition for both sets of data. In other words, there exist $t_0<t_1 < \dots <t_n$ satisfying
\begin{equation*}
\begin{array}{ll}
0\le t_0<\e \,,&\\[2mm]
t_i-t_{i-1} < \e \quad &\mbox{ for }\quad i=1,2, \dots, n \,,\\[2mm]
0\le T-t_n<\e\,, &
\end{array}
\end{equation*}
and $g_1^k, g_2^k, \dots, g_n^k \in L^2(\partial\Omega)$ such that
\begin{equation*}
\sum_{i=1}^n \int_{t_{i-1}}^{t_i} \| g^k(s)-g_i^k\|_{L^2(\partial \Omega)}\,ds < \e \,,
\end{equation*}
for $k=1,2$.
Moreover, thanks to \cite[Proposition 1.5]{BCP}, we may choose the corresponding $g_i^k=g^k(t_i)$, being each $t_i$ a Lebesgue point of $g^k$. As a consequence, $g^1(t, x)\le g^2(t, x)$ for almost all $(t,x)\in (0,T)\times \partial\Omega$ implies $g_i^1(x)\le g_i^2(x)$ for almost all $x\in \partial\Omega$ and for $i=1,\dots , n$.

Consider now the systems
\begin{equation*}
\frac{\omega_{i}^k-\omega_{i-1}^k}{t_i-t_{i-1}} + \b(\omega_i^k) \ni g_i^k \qquad\mbox{ for }\quad i=1,2,\dots,n \quad k=1,2\,,
\end{equation*}
so that
\begin{equation*}
  \omega_{i-1}^k+(t_i-t_{i-1})g_i^k\in \big(I+(t_i-t_{i-1})\b\big)(\omega_i^k)\quad\mbox{ for }\quad i=1,2,\dots,n \quad k=1,2\,.
\end{equation*}
Since $\omega_0^1(x)\le \omega_0^2(x)$ and $g_i^1(x)\le g_i^2(x)$ for almost all $x\in \partial\Omega$ and for $i=1,\dots , n$, and each resolvent $\big(I+(t_i-t_{i-1})\b\big)^{-1}$ is order preserving (see Remark \ref{comp3}), an appeal to induction leads to $\omega_i^1(x)\le \omega_i^2(x)$ for almost all $x\in \partial\Omega$ as well as $u_i^1(x)\le u_i^2(x)$ for almost every $x\in\Omega$ and for $i=1,\dots , n$.

Denoting by $\omega_\epsilon^k$ the solution to the $\epsilon$-discretization corresponding to data $g_i^k$ and $\omega_0^k$, it follows that
$\omega_\epsilon^1(t, x)\le \omega_\epsilon^2(t, x)$ for almost all $(t,x)\in (0,T)\times \partial\Omega$. Having in mind
 \begin{equation*}
   \omega_\e^k \to \omega^k \qquad \mbox{ in } L^\infty([0,T];L^2(\partial\Omega))
\end{equation*}
for $k=1,2$ (see Remark \ref{convergencia.mild-solutions}), this fact implies $\omega^1(t, x)\le \omega^2(t, x)$ for almost all $(t,x)\in (0,T)\times \partial\Omega$.
\end{pf}

\section{Existence of strong solutions}

In this Section, we are proving that the mild solution we have obtained in the previous Section is actually a strong solution to our problem. First, we introduce the concept of strong solution in our framework.

\begin{Definition}\label{def.sol.Dynamic-problem}
Let $g\in L^1(0,T;L^2(\partial\Omega))$ and let $\omega_0\in L^2(\partial\Omega)$.
We say that the pairing $(u,\omega)$ is a strong solution to problem \eqref{Dynamic-problem} if \mbox{$u \in L^\infty(0,T;BV(\Omega)) \cap L^\infty(0,T;L^2(\Omega))$} and $\omega \in C([0,T]; L^2(\partial \Omega))\cap W^{1,1}(0,T;L^2(\partial \Omega))$ such that $\omega(0)=\omega_0$ and there exists a vector field $\z \in L^\infty((0,T)\times\Omega;\R^N)$ with $\|\z\|_{L^\infty((0,T)\times\Omega)} \le 1$ satisfying  the following conditions:
\begin{equation*}
\begin{array}{ll}
(i) & \lambda u(t)-\Div(\z(t)) = 0  \qquad \mbox{ in } \dis(\Omega) \,, \\[2mm]
(ii) & (\z(t),Du(t))=|Du(t)| \qquad \mbox{ as measures in } \Omega \,, \\[2mm]
(iii) & [\z(t),\nu] = g(t)-\omega_t(t) \qquad \mbox{ for almost every } x \in \partial \Omega \,, \\[2mm]
(iv) & [\z(t),\nu]\in \sg(\omega(t)-u(t)) \qquad \mbox{ on } \partial \Omega \,,
\end{array}
\end{equation*}
for almost every $t \in (0,T)$.

Given $g\in L_{loc}^1(0,+\infty;L^2(\partial\Omega))$ and  $\omega_0\in L^2(\partial\Omega)$, we say that $(u,\omega)$  is a global strong solution to problem \eqref{Dynamic-globalproblem}
 if it is a strong solution to \eqref{Dynamic-problem} for every $T>0$.
\end{Definition}

As mentioned above, functions $u,\,\omega,\, \z,\,g$ depend on two variables: $t$ and $x$. For the sake of simplicity, most of the time we will write $u(t)$, $\omega(t)$, $\z(t)$ and $g(t)$ instead of $u(t,x)$, $\omega(t,x)$, $\z(t,x)$ and $g(t,x)$.

\begin{Theorem}
Let $\lambda >0$, and let $g \in L_{loc}^1(0,+\infty;L^2(\partial\Omega))$ and $\omega_0 \in L^2(\partial \Omega)$. Then there exists a global strong solution $(u, \omega)$ to problem \eqref{Dynamic-globalproblem}.

Furthermore, the following estimates hold:
\begin{equation}\label{est1}
  \|\omega\|_{L^\infty (0,T; L^2(\partial\Omega))}\le \|\omega_0\|_{L^2(\partial\Omega)}+\|g\|_{L^1 (0,T; L^2(\partial\Omega))}\quad\hbox{for every }T>0 \,,
\end{equation}
\begin{equation}\label{est2}
  \lambda\|u(t)\|_{L^2(\Omega)}^2+\|u(t)\|_{BV(\Omega)}
  \le 2\|\omega(t)\|_{L^1(\partial\Omega)}\quad\hbox{for almost all }t>0\,.
\end{equation}
\end{Theorem}

\begin{pf} First fix $T>0$. Applying Theorem \ref{ex-mild}, there exists a mild solution $\omega$ to the abstract Cauchy problem  $\omega_t+\b(\omega)\ni g$, $\omega(0)=\omega_0$ on $[0,T]$ with auxiliary function $u$.  We are seeing that $(u, \omega)$ is actually a strong solution.

We will divide the proof in several steps.

\medskip
\textbf{STEP 1: Solutions to $\e$-discretizations.}
\medskip

Since $\omega\in C([0,T]; L^2(\partial\Omega))$ is a mild solution, we may choose a family of \hbox{$\e$-discretizations} of $\omega_t+\b(\omega)\ni g$, $\omega(0)=\omega_0$ on $[0,T]$, in such a way that their solutions $\omega_\e$ satisfy
\begin{equation}\label{conv-omega}
  \omega_\e\to \omega \quad \mbox{ strongly in } \quad L^\infty(0,T;L^2(\partial\Omega)) \,.
\end{equation}

We will detail our notation. Fixed $0< \e \le 1$, there exists a partition $0=t_0 <t_1<\dots <t_n<T$ such that $T-t_n<\e$ and $t_i-t_{i-1} < \e$ for every $i=1,2, \dots, n$, and there exist functions $\widehat{g}_1,\dots,\widehat{g}_n \in L^2(\partial\Omega)$ such that
\begin{equation}\label{aprox-g}
\sum_{i=1}^n \int_{t_{i-1}}^{t_i} \left(\int_{\partial \Omega}|g(t,x)-\widehat{g}_i(x)|^2 \,d\h\right)^{\frac{1}{2}}\,dt <\e \,,
\end{equation}
 and so the system
 \begin{equation*}
\frac{\omega_{i}-\omega_{i-1}}{t_i-t_{i-1}} + \b(\omega_i)\ni \widehat{g}_i \quad \mbox{ for every } i=1,2,\dots,n
\end{equation*}
is an $\e$-discretization of $\omega_t+\b(\omega)\ni g$, $\omega(0)=\omega_0$ on $[0,T]$.

 We denote $\e_i = t_i-t_{i-1}$. Observe that, splitting the intervals if necessary, there is not loss of generality in assuming $\e_1 > \e_2 > \dots > \e_{n-1}$. Hence, if $t\in ]t_{i-1}, t_i]$, then $t-\e_i\in ]t_{i-2}, t_{i-1}]$.

We also define $g_\e(t,x)=\widehat{g}_i(x)$ if $t \in \, ]t_{i-1},t_i]$, for $i=1,2,\dots,n$.
Therefore, the condition \eqref{aprox-g} becomes
\begin{equation*}
\int_0^T \left(\int_{\partial \Omega}|g(t,x)-g_\e(t,x)|^2 \,d\h\right)^{\frac{1}{2}}\,dt <\e \,,
\end{equation*}
and we have the following convergence:
\begin{equation}\label{conv-g-L1enL2}
g_\e \to g \quad \mbox{ strongly in } \quad L^1(0,T;L^2(\partial\Omega)) \,.
\end{equation}

Now, the solution to the $\e$-discretization satisfies
\begin{equation*}
\omega_\e(t,x)=
\omega_i(x) \quad \mbox{if} \quad t \in \, ]t_{i-1},t_i] \quad \mbox{ for } i=1,2,\dots,n \,,
\end{equation*}
where
\begin{equation*}
\widehat{g}_i+\frac{\omega_{i-1}-\omega_i}{\e_i} \in \b(\omega_i) \quad \mbox{ for every } i=1,2,\dots,n\,.
\end{equation*}
Due to the definition of the operator $\b$, for each $i=1,2,\dots,n$, it holds
\begin{equation*}
\begin{array}{cl}
(i) &  \omega_i \in L^2(\partial\Omega) \,,\\[3mm]
(ii) &  \widehat{g}_i+\dfrac{\omega_{i-1}-\omega_i}{\e_i} \in L^\infty(\partial\Omega) \quad \mbox{ with }\quad \left\| \widehat{g}_i+\frac{\omega_{i-1}-\omega_i}{\e_i} \right\|_{L^\infty(\partial\Omega)} \le 1 \,,\\[3mm]
(iii) &  \mbox{there exists }\quad u_i \in BV(\Omega)\cap L^2(\Omega) \,,\\[3mm]
(iv) & \mbox{there exists }\quad \z_i \in L^\infty (\Omega;\R^N) \quad \mbox{ with } \quad \|\z_i\|_{L^\infty(\Omega)} \le 1 \\[3mm]
\end{array}
\end{equation*}
satisfying the following conditions
\begin{align}
\nonumber & \lambda u_i - \Div \z_i = 0 \quad \mbox{ in } \dis(\Omega) \,,\\[3mm]
\label{cond-medidas} & (\z_i,Du_i)=|Du_i| \quad \mbox{ as measures in } \Omega \,,\\[3mm]
\label{cond-signo1}  & [\z_i,\nu]= \widehat{g}_i+\dfrac{\omega_{i-1}-\omega_i}{\e_i} \quad \h\mbox{-a.e. on } \partial \Omega \,, \\[3mm]
\label{cond-signo2} & [\z_i,\nu] \in \sg(\omega_i-u_i) \quad \h\mbox{-a.e. on } \partial \Omega \,, \\[3mm]
\label{cond-distribuciones} & \displaystyle \lambda\int_\Omega u_i\,\varphi \,dx + \int_\Omega(\z_i,D\varphi) = \int_{\partial\Omega}\left(\widehat{g}_i+\dfrac{\omega_{i-1}-\omega_i}{\e_i}\right)\varphi \,d\h \,,
\end{align}
for every $\varphi \in BV(\Omega)\cap L^2(\Omega)$.

Finally, given $\e >0$, we define the following step functions:
\begin{equation*}
\begin{array}{l}
u_\e (t,x)=u_i(x) \quad \mbox{ if }\; t\in \, ]t_{i-1},t_i]\quad \mbox{ for } i=1,2,\dots,n \,,\\[3mm]
\z_\e (t,x)=\z_i(x) \quad \mbox{ if }\; t\in \,]t_{i-1},t_i]\quad \mbox{ for } i=1,2,\dots,n \,.
\end{array}
\end{equation*}

We remark that all the above step functions are defined in $[0,t_n]$. To avoid lack of definiteness, we can extend them to $]t_n,T]$ giving their value at the point $t_n$.

\medskip
\textbf{STEP 2: Existence of ${\bf \omega_t}$ in the sense of distributions.}
\medskip

Due to Definition \ref{def.operator-B} we know that
\begin{equation*}
\left\| \frac{\omega_i-\omega_{i-1}}{\e_i}-\hat{g}_i \right\|_{L^\infty(\partial\Omega)} \le 1 \qquad \mbox{ for every } i=1,2,\dots,n\,,
\end{equation*}
where $\e_i = t_i-t_{i-1}$.
Denoting $\e(t)=\e_i$ for $t\in ]t_{i-1}, t_i]$, the following equivalent bound holds:
\begin{equation}\label{acot}
\left\| \frac{\omega_\e(t)-\omega_\e(t-\e(t))}{\e(t)}-g_\e(t) \right\|_{L^\infty(\partial\Omega)} \le 1 \qquad \mbox{ for every } t \in (\e_1,T)\subset (\e, T)\,.
\end{equation}

Setting $\eta>0$, let $0<\e<\eta$ and  $t \in (\eta,T)$ be fixed. We will assume that this given $t$ satisfies
\begin{equation}\label{conv-g-L2cpp}
  g_\e(t)\to g(t)\qquad \hbox{strongly in }L^2(\partial\Omega)\,,
\end{equation}
which is a straightforward consequence of \eqref{conv-g-L1enL2}.

 Since the sequence $\left\{ \dfrac{\omega_{\e}(t)-\omega_{\e}(t-{\e(t)})}{{\e(t)}}-g_{\e}(t) \right\}$ is bounded in $L^\infty(\partial\Omega)$, there exists a subsequence and there exists a function $\rho(t) \in L^\infty(\partial\Omega)$ such that
\begin{equation}\label{conv-debil-estrella}
\frac{\omega_{\e}(t)-\omega_{\e}(t-{\e(t)})}{{\e(t)}}-g_{\e}(t)  \rightharpoonup \rho(t) \qquad *\mbox{-weakly in } L^\infty(\partial\Omega) \,.
\end{equation}
Therefore, for every $\varphi \in L^2(\partial\Omega)$ we apply \eqref{conv-debil-estrella} and \eqref{conv-g-L2cpp} to get
\begin{align*}
& \int_{\partial \Omega} \rho(t)\varphi \,d\mathcal{H}^{N-1}= \lim_{{\e} \to 0^+} \left[ \int_{\partial\Omega} \frac{\omega_{\e}(t)-\omega_{\e}(t-\e(t))}{{\e(t)}}\varphi\,d\h- \int_{\partial\Omega} g_{\e}(t)\varphi\,d\h\, \right]
\\ &\; = \lim_{{\e} \to 0^+}  \int_{\partial\Omega} \frac{\omega_{\e}(t)-\omega_{\e}(t-{\e(t)})}{{\e(t)}}\varphi\,d\h- \int_{\partial\Omega} g(t)\varphi\,d\h \,.
\end{align*}
Then, we have
\begin{equation*}
\lim_{{\e} \to 0^+} \int_{\partial\Omega} \frac{\omega_{\e}(t)-\omega_{\e}(t-{\e(t)})}{{\e(t)}}\varphi\,d\h= \int_{\partial\Omega} (g(t)+\rho(t))\varphi\,d\h \,,
\end{equation*}
that is,
\begin{equation*}
\frac{\omega_{\e}(t)-\omega_{\e}(t-{\e(t)})}{{\e(t)}} \rightharpoonup g(t)+\rho(t) \qquad \mbox{weakly in } L^2(\partial\Omega)\,.
\end{equation*}
\\
We take now the function $\psi \in C_0^1(0,T;L^2(\partial \Omega))$ such that $\sop \psi \subseteq \,]\e,T-\e[$, obtaining
\begin{align*}
& \int_0^T  \int_{\partial \Omega} \frac{\omega_\e(t)-\omega_\e(t-\e(t))}{\e(t)}\psi(t) \,d\h\, dt
\\ &\; = \int_0^T\int_{\partial \Omega} \frac{\omega_\e(t)}{\e(t)}\psi(t) \,d\h\, dt -  \int_0^T\int_{\partial \Omega} \frac{\omega_\e(t-\e(t))}{\e(t)}\psi(t)\,d\h\, dt
\\ &\; = \int_0^T\int_{\partial \Omega} \frac{\omega_\e(t)}{\e(t)}\psi(t)\,d\h\, dt-  \int_0^T\int_{\partial \Omega} \omega_\e(t)\frac{\psi(t+\e(t))}{\e(t)}\,d\h\, dt
\\ &\; = - \int_0^T\int_{\partial \Omega}\omega_\e(t) \frac{\psi(t+\e(t))-\psi(t)}{\e(t)}\,d\h\, dt \,.
\end{align*}
On the other hand, having in mind \eqref{conv-g-L1enL2} and \eqref{conv-debil-estrella} (and also \eqref{acot}), it follows that
\begin{align*}
& \int_0^T\int_{\partial \Omega} (g(t)+\rho(t))\psi(t) \,d\h\, dt
=\lim_{\e \to 0^+} \int_0^T\int_{\partial\Omega} g_\e(t)\psi(t)\,d\h\, dt
\\ &\; + \lim_{\e \to 0^+}\int_0^T\int_{\partial\Omega} \Big(\frac{\omega_\e(t)-\omega_\e(t-\e(t))}{\e(t)}-g_\e(t)\Big)\psi(t)\,d\h\, dt
\\ &\; =\lim_{\e \to 0^+} \int_0^T\int_{\partial\Omega} \frac{\omega_\e(t)-\omega_\e(t-\e(t))}{\e(t)}\psi(t)\,d\h\, dt\,.
\end{align*}
Therefore,
\begin{align*}
& \int_0^T \int_{\partial \Omega} (g(t)+\rho(t))\psi(t) \,d\h\, dt
\\ &\; =\lim_{\e \to 0^+} \int_0^T\int_{\partial\Omega} \frac{\omega_\e(t)-\omega_\e(t-\e(t))}{\e(t)}\psi(t)\,d\h\, dt
\\ &\; =- \lim_{\e \to 0^+} \int_0^T\int_{\partial \Omega}\omega_\e(t) \frac{\psi(t+\e(t))-\psi(t)}{\e(t)}\,d\h\, dt
\\ &\; = -\int_0^T\int_{\partial \Omega} \omega(t)\psi_t(t)\,d\h\,dt \,,
\end{align*}
due to \eqref{conv-omega}.
Then, the distributional derivative of $\omega$ is $\omega_t=g+\rho \in L^1(0,T; L^2(\partial\Omega))$ and  it also holds $\|\omega_t(t)-g(t)\|_{L^\infty(\partial\Omega)} \le 1$ for almost every $t\in (0,T)$.
\\
Moreover, we have that
\begin{equation}\label{conv-omega_t}
\frac{\omega_\e(t)-\omega_\e(t-\e(t))}{\e(t)} \rightharpoonup g(t)+\rho(t)=\omega_t(t) \qquad \mbox{weakly in } L^2(\partial\Omega)\,.
\end{equation}

We point out that, since the operator $\b$ is $m$-accretive, function $\omega$ is absolutely continuous and differentiable in almost every $t \in (0,T)$ and besides it is a mild solution to problem $\omega_t +\b(\omega) \ni g$ on $(0,T)$, it yields that function $\omega$ is also a strong solution (see \cite[Theorem 7.1]{BCP}). In other words, $g(t)-\omega_t(t) \in \b(\omega(t))$ holds for almost every $t \in (0,T)$. This concludes the proof in what the boundary concerns, which is where the semigroup is defined. Hence,  for every $t\in(0,T)$ fixed, there exist an auxiliary $BV$-function and a vector field satisfying Definition \ref{def.operator-B}.
Nevertheless, in the domain $(0,T)$ there may be a problem of measurability since the strong solution only provide us the functions pointwise in time. In the sequel, we will find $u\in L^\infty(0,T;BV(\Omega))\cap L^\infty(0,T;L^2(\Omega))$ and $\z\in L^\infty((0,T)\times\Omega; \R^N)$ satisfying all the requirements of Definition \ref{def.sol.Dynamic-problem}.

\medskip
\textbf{STEP 3: Existence of $\z\in L^\infty((0,T)\times\Omega; \R^N)$.}
\medskip

This fact is an easy consequence of $\|\z_\e\|_{L^\infty((0,T)\times\Omega)}\le1$ for all $\e>0$. Then there exists a vector field $\z\in L^\infty((0,T)\times\Omega;\R^N)$ such that, up to subsequences,
\begin{equation}\label{conv-campo}
\z_\e \tod \z \qquad *\mbox{-weakly in } L^\infty((0,T)\times\Omega)\,.
\end{equation}

\vskip 3cm
\textbf{STEP 4: The sequence $\{{\bf u_\e}\}$ is bounded  in ${\bf L^2(0,T;L^2(\Omega))}$ and  in ${\bf L^1(0,T;BV(\Omega))}$.}
\medskip

We take $u_i \in BV(\Omega)\cap L^2(\Omega)$ as test function in \eqref{cond-distribuciones} and we get
\begin{align*}
\lambda & \int_\Omega u_i^2 \,dx + \int_\Omega (\z_i,Du_i) = \int_{\partial\Omega}\left(\widehat{g}_i+\dfrac{\omega_{i-1}-\omega_i}{\e_i}\right) u_i \,d\h
\\ & =\int_{\partial \Omega} \left(\widehat{g}_i+\dfrac{\omega_{i-1}-\omega_i}{\e_i}\right)(u_i-\omega_i) \,d\h + \int_{\partial\Omega} \left(\widehat{g}_i+\dfrac{\omega_{i-1}-\omega_i}{\e_i}\right) \omega_i \,d\h
\\ & =-\int_{\partial \Omega} |u_i-\omega_i| \,d\h + \int_{\partial\Omega} \widehat{g}_i \omega_i \,d\h + \int_{\partial\Omega}\dfrac{\omega_{i-1}-\omega_i}{\e_i} \omega_i \,d\h\,,
\end{align*}
where we have used conditions \eqref{cond-signo1} and \eqref{cond-signo2}. Then,  condition \eqref{cond-medidas} implies
\begin{multline}\label{ec-1}
\lambda  \int_\Omega u_i^2 \,dx + \int_\Omega |Du_i| +\int_{\partial \Omega} |u_i-\omega_i| \,d\h+\int_{\partial\Omega} \frac{\omega_i-\omega_{i-1}}{\e_i}\omega_i \,d\h
\\  =\int_{\partial\Omega} \widehat{g}_i \omega_i \,d\h \,,
\end{multline}
and, dropping nonnegative terms, we get the following inequality
\begin{equation}\label{ec-1.1}
\int_{\partial\Omega} \frac{\omega_i-\omega_{i-1}}{\e_i}\omega_i \,d\h\le \int_{\partial\Omega} \widehat{g}_i \, \omega_i \,d\h \,,
\end{equation}
for every $i=1,2,\dots,n$.

Next, we will show that $\displaystyle \int_{\partial\Omega} \omega_i^2 \,d\h$ is bounded by a constant which does not depend on $i$.
Using H\"older's inequality, condition \eqref{ec-1.1} and then H\"older's inequality again, we get
\begin{align*}
& \int_{\partial\Omega} \omega_i^2 \,d\h - \left( \int_{\partial\Omega}\omega_i^2 \,d\h\right)^\frac{1}{2}\left(\int_{\partial\Omega}\omega^2_{i-1}\,d\h\right)^\frac{1}{2}
\\ &\; \le \int_{\partial\Omega}\omega^2_i \,d\h - \int_{\partial\Omega}\omega_i\,\omega_{i-1} \,d\h = \int_{\partial\Omega} (\omega_i-\omega_{i-1})\,\omega_i \,d\h \\ &\; \le \e_i \int_{\partial\Omega} \widehat{g}_i\omega_i \,d\h \le \e_i \left(\int_{\partial\Omega}\widehat{g}_i^2\,d\h\right)^\frac{1}{2}\left(\int_{\partial\Omega}\omega_i^2\,d\h\right)^\frac{1}{2} \,.
\end{align*}
Now, if $\displaystyle \int_{\partial\Omega}\omega_i^2 \,d\h \not=0$, we divide the previous inequality by $\displaystyle \left(\int_{\partial\Omega}\omega_i^2 \,d\h\right)^\frac{1}{2} $ and so we get
\begin{equation*}
\left(\int_{\partial\Omega}\omega_i^2\,d\h\right)^\frac{1}{2}-\left(\int_{\partial\Omega} \omega^2_{i-1} \,d\h\right)^\frac{1}{2} \le \e_i \left(\int_{\partial\Omega} \widehat{g}_i^2\,d\h\right)^\frac{1}{2}\,,
\end{equation*}
for every $i=1,2,\dots,n$.
We fix now $i\in \{ 1, 2,\dots, n \}$ and sum the previous inequality for $k=1,2, \dots,i-1,i$:
\begin{align*}
& \bigg(\int_{\partial\Omega} \omega_i^2\,d\h\bigg)^\frac{1}{2}-\left(\int_{\partial\Omega}\omega_0^2\,d\h\right)^\frac{1}{2}
\\ &\; =\sum_{k=1}^i\left[\left(\int_{\partial\Omega}\omega_k^2\,d\h\right)^\frac{1}{2}-\left(\int_{\partial\Omega} \omega^2_{k-1} \,d\h\right)^\frac{1}{2}\right]
\\ &\; \le \sum_{k=1}^i \e_k \left(\int_{\partial\Omega} \widehat{g}_k^2\,d\h\right)^\frac{1}{2}=\sum_{k=1}^i \int_{t_{k-1}}^{t_k} \left(\int_{\partial\Omega} \widehat{g}_k^2\,d\h\right)^\frac{1}{2} dt\,.
\end{align*}
We can perform easy manipulations to get
\begin{align}\label{ec-2.1}
& \bigg(\int_{\partial\Omega} \omega_i^2\,d\h\bigg)^\frac{1}{2} \le \left(\int_{\partial\Omega}\omega_0^2\,d\h\right)^\frac{1}{2} + \int_{0}^{t_i}\left(\int_{\partial\Omega} g_\e^2\,d\h\right)^\frac{1}{2}dt
\\\nonumber &\; \le \left(\int_{\partial\Omega}\omega_0^2\,d\h\right)^\frac{1}{2} + \int_{0}^T\left(\int_{\partial\Omega} g_\e^2\,d\h\right)^\frac{1}{2}dt \,.
\end{align}
Therefore, we deduce that
\begin{align*}
\|\omega_i\|_{L^2(\partial\Omega)} & \le \|\omega_0\|_{L^2(\partial\Omega)} + \|g_\e\|_{L^1(0,T;L^2(\partial\Omega))}
\\ & \le \|\omega_0\|_{L^2(\partial\Omega)} + 1 +\|g\|_{L^1(0,T;L^2(\partial\Omega))}= M \,,
\end{align*}
 for every  $i=1,2,\dots,n $.
That is, the sequence $\omega_\e(t,x)$ is bounded in $L^\infty(0,T;L^2(\partial\Omega))$.
In addition, since $g_\e(t),\,\omega_\e(t) \in L^2(\partial\Omega)$, we can use H\"older's inequality to get
\begin{align}\label{ec-3}
& \int_0^T \int_{\partial\Omega} g_\e\omega_\e\,d\h \,dt\le \int_0^T \left(\int_{\partial\Omega} g_\e^2\,d\h\right)^\frac{1}{2}\left(\int_{\partial\Omega} \omega_\e^2\,d\h\right)^\frac{1}{2} \,dt
\\ &\; \le  \int_0^T M \left(\int_{\partial\Omega} g_\e^2\,d\h\right)^\frac{1}{2} \,dt =M \int_0^T \left(\int_{\partial\Omega} g_\e^2\,d\h\right)^\frac{1}{2} \,dt\le M^2 \,.
\end{align}

On the other hand, we know that $\dfrac{\omega_i^2-\omega^2_{i-1}}{2}\le (\omega_i-\omega_{i-1})\,\omega_i$ and we deduce from \eqref{ec-1} that
\begin{align*}
\lambda & \int_\Omega u_i^2 \,dx + \int_\Omega |Du_i| +\int_{\partial \Omega} |u_i-\omega_i| \,d\h
\\ & +\frac{1}{\e_i}\int_{\partial\Omega} \frac{\omega_i^2}{2} \,d\h-\frac{1}{\e_i}\int_{\partial\Omega}\frac{\omega^2_{i-1}}{2}\,dx
 \le\int_{\partial\Omega} \widehat{g}_i \omega_i \,d\h \,.
\end{align*}
Now, we integrate the previous inequality between $t_{i-1}$ and $t_i$ to obtain
\begin{align*}
\lambda & \int_{t_{i-1}}^{t_i}\int_\Omega u_i^2 \,dx\,dt +\int_{t_{i-1}}^{t_i}\int_\Omega |Du_i|\,dt +\int_{t_{i-1}}^{t_i}\int_{\partial \Omega} |u_i-\omega_i| \,d\h\,dt
\\ & +\int_{\partial\Omega} \frac{\omega_i^2}{2} \,d\h-\int_{\partial\Omega}\frac{\omega^2_{i-1}}{2}\,dx
 \le \int_{t_{i-1}}^{t_i}\int_{\partial\Omega} \widehat{g}_i \omega_i \,d\h \,dt\,,
\end{align*}
for every $i=1,2,\dots,n$. Adding all these terms from $i=0$ to $n$:
\begin{align*}
\lambda & \sum_{i=1}^n\int_{t_{i-1}}^{t_i}\int_\Omega u_i^2 \,dx\,dt +\sum_{i=1}^n\int_{t_{i-1}}^{t_i}\int_\Omega |Du_i|\,dt +\sum_{i=1}^n\int_{t_{i-1}}^{t_i}\int_{\partial \Omega} |u_i-\omega_i| \,d\h\,dt
\\ & +\int_{\partial\Omega} \frac{\omega_n^2}{2} \,d\h-\int_{\partial\Omega}\frac{\omega^2_0}{2}\,dx
 \le \sum_{i=1}^n\int_{t_{i-1}}^{t_i}\int_{\partial\Omega} \widehat{g}_i \omega_i \,d\h\,dt\,.
\end{align*}
Then, we finally get
\begin{align*}
\lambda & \int_0^T\int_\Omega u_\e^2 \,dx\,dt +\int_0^T \int_\Omega |Du_\e|\,dt +\int_0^T\int_{\partial \Omega} |u_\e-\omega_\e| \,d\h\,dt
+\int_{\partial\Omega} \frac{\omega_n^2}{2} \,d\h
\\ &  \le \int_{\partial\Omega}\frac{\omega^2_0}{2}\,dx
+ \int_0^T\int_{\partial\Omega} g_\e \omega_\e \,d\h \,dt \le 2M^2 \,,
\end{align*}
where we have used \eqref{ec-3}.

Therefore, we have proved that
\begin{equation*}
\lambda\int_0^T\int_\Omega u_\e^2 \,dx\,dt +\int_0^T\int_\Omega |Du_\e|\,dt \le 2M^2 \,.
\end{equation*}
That is, the sequence $\{ u_\e \}$ is bounded in $L^2(0,T;L^2(\Omega))$ and, by H\"older's inequality, it is also bounded in $L^1(0,T;BV(\Omega))$.

As a first consequence, there exists a measurable function $u \in L^2((0,T)\times\Omega)$ such that
\begin{equation}\label{conv-sol}
u_\e \rightharpoonup u \qquad\mbox{ weakly in } \quad L^2((0,T)\times\Omega) \,.
\end{equation}

\medskip
\textbf{STEP 5: Function ${\bf u}$ belongs to ${\bf L^\infty (0,T;BV(\Omega))
\cap L^\infty(0,T;L^2(\Omega))}$.}
\medskip

Let $t\in (0,T)$. Observe that \eqref{ec-1}, written in terms of the
approximate solutions, becomes
\begin{align*}
\lambda & \int_\Omega u_\e(t)^2 \,dx + \int_\Omega |Du_\e(t)| +\int_{\partial
\Omega} |u_\e(t)-\omega_\e(t)| \,d\h
\\ & =\int_{\partial\Omega} \left( g_\e(t)-\frac{\omega_\e(t)-\omega_\e(t-\e_i)
}{\e_i}\right) \omega_\e(t) \,d\h  \le \int_{\partial\Omega} |\omega_\e(t)| \,d\h
 \,,
\end{align*}
because of \eqref{cond-signo1}.
It follows that
\begin{equation}\label{cotaBV}
\lambda\int_\Omega u_\e(t)^2 \,dx + \int_\Omega |Du_\e(t)| \le M_1=M\mathcal H^{N-1}(\partial\Omega)^{1/2} \,,
\end{equation}
so that the sequence $\{ u_\e\}$ is bounded in $L^\infty(0,T;L^2(\Omega))$ and, on account of H\"older's inequality, in $L^\infty(0,T;BV(\Omega))$.

In order to see that $u\in L^\infty(0,T;BV(\Omega))\cap L^\infty(0,T;L^2(
\Omega))$, we need to let $\e\to0$ in the above inequality. To this end, we first fix
$\xi\in C_0^\infty((0,T))$ such that $\xi\ge 0$ and observe that, for each $v\in L^\infty(0,T; BV(\Omega))$,
\begin{equation*}
\int_0^T\int_\Omega |Dv(t,x)| \xi(t)\, dx\,dt=\sup\left\{\int_0^T\int_\Omega v(t,x)\, \Div \psi(x)\,\xi(t)\, dx \, dt\right\}
\end{equation*}
where the supremum is taken among all $\psi\in C_0^1(\Omega;\R^N)$ such that
$|\psi(x)|\le 1$. Since for every $\xi$ and $\psi$ fixed we have the continuity
of
\begin{equation*}
v\mapsto \int_0^T\int_\Omega v(t,x)\, \Div \psi(x)\,\xi(t)\, dx \, dt
\end{equation*}
with respect to the weak convergence in $L^1((0,T)\times \Omega)$, we deduce
that the functional
\begin{equation*}
v\mapsto \int_0^T\int_\Omega |Dv(t,x)| \xi(t)\, dx\, dt
\end{equation*}
is lower semicontinuous with respect to the weak convergence in $L^1((0,T)
\times \Omega)$. Thus, $u_\e\rightharpoonup u$ weakly in $L^2((0,T)\times \Omega)$
implies
\begin{equation*}
\int_0^T\int_\Omega |Du(t,x)| \xi(t)\, dx\, dt\le\liminf_{\e\to
0}\int_0^T\int_\Omega |Du_\e(t,x)| \xi(t)\, dx\, dt\,.
\end{equation*}
It follows that $u\in L^\infty(0,T;BV(\Omega))\cap L^\infty(0,T;L^2(\Omega))$
and
\begin{equation*}
\lambda\int_\Omega u(t)^2 \,dx + \int_\Omega |Du(t)| \le M_1 \quad\hbox{for
almost all }t\in(0,T)\,.
\end{equation*}

\medskip
\textbf{STEP 6: ${\bf \lambda u(t)-\Div \z(t) = 0}$ holds in $\mathcal D'(\Omega)$ for almost every ${\bf t \in (0,T)}$.}
\medskip

Observe that, since \eqref{cond-distribuciones} holds for every  $i=1,2,\dots,n$,  we have
\begin{equation*}
\lambda\int_\Omega u_\e(t,x)\,\varphi(x) \,dx + \int_\Omega\z_\e(t,x)\cdot\nabla\varphi(x)\, dx = 0 \,,
\end{equation*}
for every $\varphi \in C_0^\infty(\Omega)$ and for every $t\in (0,T)$. Considering $\xi\in C_0^\infty(0,T)$, it follows from convergences \eqref{conv-sol} and \eqref{conv-campo} that
\begin{equation*}
\lambda\int_0^T\int_\Omega u(t, x)\,\varphi(x)\, \xi(t) \,dx \, dt + \int_0^T\int_\Omega\z(t,x)\cdot\nabla\varphi(x)\, \xi(t)\, dx\,  dt = 0 \,.
\end{equation*}
Therefore, for almost every $t\in (0,T)$, we deduce
\begin{equation*}
\lambda\int_\Omega u(t)\,\varphi \,dx  + \int_\Omega\z(t)\cdot\nabla\varphi\, dx= 0 \,,
\end{equation*}
and so Step 6 is proved and $\Div\z(t)\in L^2(\Omega)$ for almost all $t\in (0,T)$.

\medskip
\textbf{STEP 7: Proof of ${\bf \omega_t(t)+[\z(t), \nu]=g(t)}$ on $\partial\Omega$ holds for almost all ${\bf t\in (0,T)}$.}
\medskip

As a consequence of $\Div\z(t)\in L^2(\Omega)$, we may apply Green's formula to the vector field $\z(t)$. So Step 3 implies the following limit
\begin{equation}\label{conv-campo2}
[\z_\e(t,x),\nu(x)] \tod  [\z(t,x),\nu(x)] \qquad *\mbox{-weakly in } L^\infty((0,T)\times\partial\Omega)\,.
\end{equation}
On the other hand,
\begin{equation}\label{conv-traza}
-[\z_\e(t), \nu]=\frac{\omega_\e(t)-\omega_\e(t-\e(t))}{\e(t)}-g_\e(t)\rightharpoonup \rho(t) \qquad *\mbox{-weakly in } L^\infty(\partial\Omega) \,.
\end{equation}
Taking $\varphi\in L^2(\partial\Omega)$ and $\xi\in L^2((0,T))$, we may compute the limit of
\begin{equation*}
  \int_0^T\int_{\partial\Omega} [\z_\e(t,x), \nu(x)]\varphi(x)\xi(t)\, d\mathcal H^{N-1}\, dt
\end{equation*}
 using \eqref{conv-campo2} and \eqref{conv-traza}; then
 \begin{equation*}
  \int_0^T\int_{\partial\Omega} [\z(t,x), \nu(x)]\varphi(x)\xi(t)\, d\mathcal H^{N-1}\, dt=-\int_0^T\int_{\partial\Omega} \rho(t,x)\varphi(x)\xi(t)\, d\mathcal H^{N-1}\, dt\,.
\end{equation*}
Thus, $[\z(t), \nu]=-\rho(t)$ $\mathcal H^{N-1}$-a.e. on $\partial\Omega$  and
\begin{equation}\label{conv-traza2}
[\z_\e(t), \nu]\rightharpoonup [\z(t), \nu] \qquad *\mbox{-weakly in } L^\infty(\partial\Omega)
\end{equation}
for almost all $t\in(0,T)$.
Recalling \eqref{conv-omega_t}, we also deduce that the identity
\begin{equation*}
  \omega_t(t)+[\z(t), \nu]=g(t)
\end{equation*}
holds on $\partial\Omega$ for almost all $t\in (0,T)$.

\medskip
\textbf{STEP 8: For almost every ${\bf t \in(0,T)}$ there exists a subsequence satisfying some useful convergences.}
\medskip

Let $t\in (0,T)$. From \eqref{cotaBV} it follows that
\begin{equation*}
\lambda \int_\Omega u_\e(t)^2 \,dx \le M_1 \,.
\end{equation*}
Then, $\{u_\e(t)\}$ is bounded in $L^2(\Omega)$ and there exist $\widehat u(t)\in L^2(\Omega)$ and a subsequence $\{u_{\e^t}(t)\}$ (we remark that the subsequence we find depends on $t$)  such that
\begin{equation}\label{conv-sol-v}
u_{\e^t}(t) \rightharpoonup \widehat u(t) \qquad \mbox{weakly in } L^2(\Omega) \,.
\end{equation}

Now we go back to \eqref{cotaBV} which is an estimate of $\{ u_{\e^t} (t)\}$ in $BV(\Omega)$ for a fixed $t\in (0,T)$.
Thus, there exists a further subsequence (not relabeled) such that converges to a BV--function strongly in $L^1(\Omega)$.
Since we have proved \eqref{conv-sol-v}, we conclude that
\begin{equation}\label{conv-fuerte-L1}
u_{\e^t}(t) \to \widehat u(t) \qquad \mbox{strongly in } L^1(\Omega)\,.
\end{equation}

On the other hand, fixed $t\in(0,T)$, the sequence $\{\z_{\e^t}(t)\}$ is bounded in $L^\infty(\Omega)$ since $\|\z_{\e^t}(t)\|_{L^\infty(\Omega)} \le 1$, then (passing to a subsequence if necessary) there exists a vector field $\widehat{\z}(t) \in L^\infty(\Omega)$ such that
\begin{equation}\label{conv-campo-v}
\z_{\e^t}(t) \rightharpoonup \widehat{\z}(t) \qquad *\mbox{-weakly in } L^\infty (\Omega) \,.
\end{equation}

\medskip
\textbf{STEP 9: ${\bf \lambda \widehat u(t)-\Div \widehat{\z}(t) = 0}$ holds for almost every ${\bf t \in (0,T)}$.}
\medskip

Observe that, since \eqref{cond-distribuciones} holds for every  $i=1,2,\dots,n$,  we have
\begin{equation*}
\lambda\int_\Omega u_{\e^t}(t)\,\varphi \,dx + \int_\Omega\z_{\e^t}(t)\cdot\nabla\varphi\, dx = 0 \,,
\end{equation*}
for every $\varphi \in C_0^\infty(\Omega)$ and for almost every $t\in (0,T)$. Then, it follows from \eqref{conv-sol-v} and \eqref{conv-campo-v} that
\begin{equation*}
\lambda\int_\Omega \widehat u(t)\,\varphi \,dx + \int_\Omega\widehat{\z}(t)\cdot\nabla\varphi\,dx = 0 \,,
\end{equation*}
and so $\Div\widehat{\z}(t)\in L^2(\Omega)$ for almost all $t\in (0,T)$.

We point out that, as consequence of \eqref{conv-campo-v} and Green's formula, we also get
\begin{equation}\label{conv-campo2-v}
[\z_{\e^t}(t),\nu] \rightharpoonup [\widehat{\z}(t),\nu] \qquad *\mbox{-weakly in } L^\infty (\partial\Omega) \,.
\end{equation}
Having in mind \eqref{conv-traza2}, we conclude that $[\z(t),\nu]=[\widehat{\z}(t),\nu]$ on $\partial\Omega$.

\medskip
\textbf{STEP 10:  ${\bf (\widehat{\z}(t),D\widehat u(t))=|D\widehat u(t)|}$ as measures in ${\bf \Omega}$ for almost every ${\bf t \in(0,T)}$.}
\medskip

Fix $t\in (0,T)$ such that $u_{\e^t}(t) \to \widehat u(t)$ strongly in $L^1(\Omega)$, $\z_{\e^t}(t) \to \widehat{\z}(t)$ $*$-weakly in $L^\infty(\Omega)$ and the distributional equation $\lambda\widehat u(t)-\Div\widehat{\z}(t)=0$ holds.
Given $\varphi \in C_0^\infty (\Omega)$ with $\varphi \ge 0$, we take the test function $\varphi \, u_{\e^t}(t)$ in \eqref{cond-distribuciones} and we obtain
\begin{equation}\label{ec-5}
\lambda\int_\Omega u_{\e^t}(t)^2 \varphi \,dx + \int_\Omega u_{\e^t}(t)\, \z_{\e^t}(t) \cdot \nabla\varphi \,dx+\int_\Omega\varphi \, |Du_{\e^t}(t)| =0 \,.
\end{equation}
We want to take limits when $\e^t$ goes to $0^+$ in each term of \eqref{ec-5}.
\\
On the one hand, the lower semicontinuity of the total variation (see Theorem \ref{teo.semicont-inferior}) provide us of
\begin{equation*}
\int_\Omega \varphi \,|D\widehat u(t)| \le \lim_{\e^t \to 0^+} \int_\Omega\varphi \, |Du_{\e^t}(t)| \,.
\end{equation*}
On the other hand, \eqref{conv-sol-v} implies
\begin{equation*}
\lambda\int_\Omega \widehat u(t)^2  \varphi \,dx \le \liminf_{\e^t \to 0^+} \lambda\int_\Omega u_{\e^t}(t)^2 \varphi \,dx \,.
\end{equation*}
Moreover,
\begin{equation*}
\lim_{\e^t \to 0^+} \int_\Omega u_{\e^t}(t)\, \z_{\e^t}(t) \cdot \nabla\varphi \,dx= \int_\Omega \widehat u(t)\, \widehat{\z}(t) \cdot \nabla\varphi \,dx \,.
\end{equation*}
Therefore, letting $\e^t \to 0^+$ in \eqref{ec-5} we get
\begin{equation*}
\lambda\int_\Omega \widehat u(t)^2 \varphi \,dx + \int_\Omega \widehat u(t)\, \widehat{\z}(t) \cdot \nabla\varphi \,dx+\int_\Omega\varphi \, |D\widehat u(t)| \le 0 \,,
\end{equation*}
which, using the previous step, can be written as
\begin{align*}
& \int_\Omega \varphi \, |D\widehat u(t)| \le -\lambda\int_\Omega \widehat u(t)^2 \varphi \,dx -\int_\Omega \widehat u(t)\, \widehat{\z}(t) \cdot \nabla\varphi \,dx
\\ &\; =-\int_\Omega \widehat u(t)\, \varphi \,\Div\widehat{\z}(t) -\int_\Omega \widehat u(t)\, \widehat{\z}(t) \cdot \nabla\varphi \,dx = \langle (\widehat{\z}(t),D\widehat u(t)),\varphi \rangle \,.
\end{align*}

Since this inequality holds for every $\varphi\ge0$, we have that $|D\widehat u(t)| \le (\widehat{\z}(t),D\widehat u(t))$ as measures. The reverse inequality is straightforward, so that the equality holds and Step 10 is proved.

\medskip
\textbf{STEP 11: Boundary condition: ${\bf [\z(t),\nu] \in \sgb (\omega(t)-\widehat u(t))}$ holds on $\partial\Omega$ for almost every ${\bf t \in (0,T)}$.}
\medskip

As in Step 10, fix $t\in (0,T)$ such that the previous Steps hold true and take  $u_{\e^t}(t)$ as a test function in \eqref{cond-distribuciones}; then
\begin{equation*}
\lambda\int_\Omega u_{\e^t}(t)^2 \,dx + \int_\Omega (\z_{\e^t}(t),Du_{\e^t}(t)) =\int_{\partial\Omega}u_{\e^t}(t)[\z_{\e^t}(t),\nu]\, d\h \,.
\end{equation*}
Applying \eqref{cond-medidas} and \eqref{cond-signo2}, we have
\begin{align*}
\lambda & \int_\Omega u_{\e^t}(t)^2 \,dx + \int_\Omega |Du_{\e^t}(t)|+\int_{\partial\Omega}|u_{\e^t}(t)-\omega_{\e^t}(t)|\, d\h
\\ & =\int_{\partial\Omega}\omega_{\e^t}(t)[\z_{\e^t}(t),\nu]\, d\h  \,,
\end{align*}
which leads to
\begin{align*}
\lambda & \int_\Omega u_{\e^t}(t)^2 \,dx + \int_\Omega |Du_{\e^t}(t)|+\int_{\partial\Omega}|u_{\e^t}(t)-\omega(t)|\, d\h
\\ & \le \int_{\partial\Omega}|\omega_{\e^t}(t)-\omega(t)|\, d\h+\int_{\partial\Omega}\omega_{\e^t}(t)[\z_{\e^t}(t),\nu]\, d\h  \,.
\end{align*}
To let ${\e^t}\to0$, in the first term we use \eqref{conv-sol-v}, while in the second and third terms we apply  Theorem \ref{teo.semicont-inferior}. The right-hand side is a consequence of the convergence
$\omega_{\e^t}(t)\to \omega(t)$ strongly in $L^2(\partial\Omega)$ and $[\z_{\e^t}(t),\nu]\rightharpoonup[\z(t),\nu]$ $*$-weakly in $L^\infty(\partial\Omega)$. Hence,
\begin{align}\label{ec-fron1}
\lambda & \int_\Omega \widehat u(t)^2 \,dx + \int_\Omega |D\widehat u(t)|+\int_{\partial\Omega}|\widehat u(t)-\omega(t)|\, d\h
\\ \nonumber & \le\int_{\partial\Omega}\omega(t)[\z(t),\nu]\, d\h \,.
\end{align}

On the other hand, Step 9, Step 10 and Green's formula imply
\begin{equation}\label{ec-fron2}
\lambda\int_\Omega \widehat u(t)^2 \,dx + \int_\Omega |D\widehat u(t)|=\int_{\partial\Omega}\widehat u(t)[\z(t), \nu]\, d\h\,.
\end{equation}
Combining \eqref{ec-fron1} and \eqref{ec-fron2}, it yields
\begin{equation*}
\int_{\partial\Omega}|\widehat u(t)-\omega(t)|\, d\h+\int_{\partial\Omega}(\widehat u(t)-\omega(t))[\z(t), \nu]\, d\h\le0\,,
\end{equation*}
from where Step 11 follows.

\medskip
\textbf{STEP 12: ${\bf u_\e(t)\tod \widehat u(t)}$ in $L^2(\Omega)$ and ${\bf u_\e(t) \to \widehat{u}(t)}$ in $L^1(\Omega)$ for almost every ${\bf t \in(0,T)}$.}
\medskip

Fix $t \in (0,T)$ such that the previous steps hold. We have proved that there exists a subsequence $\{u_{\e^t}\}$ and a function $\widehat u(t)\in L^2(\Omega)\cap BV(\Omega)$ such that \eqref{conv-sol-v} and \eqref{conv-fuerte-L1} hold, and $\widehat u(t)$ is a solution to the Dirichlet problem
\begin{equation*}
\left\{\begin{array}{ll}
\lambda v-\Div\left(\frac{Dv}{|Dv|}\right)=0\,, &\hbox{in }\Omega\,;\\
v=\omega(t)\,, &\hbox{on }\partial\Omega\,.
\end{array}\right.
\end{equation*}
The uniqueness of solutions to this problem implies that the whole sequence $\{u_\e(t)\}$ converges to $\widehat{u}(t)$ weakly in $L^2(\Omega)$ and strongly in $L^1(\Omega)$. We remark that we may also assume that $\{u_\e(t)\}$ converges to $\widehat{u}(t)$ a.e. in $\Omega$.

\medskip
\textbf{STEP 13: ${\bf u(t)=\widehat u(t)}$ for almost every ${\bf t \in(0,T)}$.}
\medskip

Since $u_\e$ are measurable functions in $(0,T)\times\Omega$, the pointwise limit function $\widehat{u}$ is also measurable in $(0,T)\times\Omega$.

Considering now $\varphi\in L^2(\Omega)$ and $\xi \in L^2((0,T))$, the following inequality holds
\begin{equation*}
\left|\int_\Omega u_\e(t,x)\varphi(x)\xi(t)\,dx\right|\le |\xi(t)|\left(\int_\Omega u_\e(t,x)^2 \,dx\right)^{\frac{1}{2}}\left(\int_\Omega \varphi(x)^2\,dx\right)^{\frac{1}{2}} \le K|\xi(t)|
\end{equation*}
for certain constant $K>0$, by \eqref{cotaBV}.
This inequality allows us to use the dominated convergence Theorem and obtain
\begin{equation*}
\lim_{\e\to0^+} \int_0^T \int_\Omega u_\e(t,x)\varphi(x)\xi(t)\,dx\,dt = \int_0^T \left[\lim_{\e\to0^+} \int_\Omega u_\e(t,x)\varphi(x)\xi(t)\,dx\right]\,dt\,,
\end{equation*}
so that
\begin{equation*}
\int_0^T \int_\Omega u(t,x)\varphi(x)\xi(t)\,dx\,dt = \int_0^T\int_\Omega \widehat{u}(t,x)\varphi(x)\xi(t)\,dx\,dt\,.
\end{equation*}
Therefore, we get that $u(t,x)=\widehat{u}(t,x)$ for almost every $(t,x)\in(0,T)\times\Omega$.

\medskip
\textbf{STEP 14: The pairing $(u,\omega)$ is a strong solution to problem \eqref{Dynamic-problem}.}
\medskip

Having in mind Steps 6, 7, 11 and 13, it only remains to check the equality $(\z(t), Du(t))=|Du(t)|$ for almost all $t\in (0,T)$. Now, a remark is in order. By Steps 10 and 13, we already know that $(\widehat \z(t), Du(t))=|Du(t)|$ holds for almost all $t\in (0,T)$. Nevertheless, the way we have obtained the vector field $\widehat \z$ does not imply that it is measurable in $(0,T)\times\Omega$. Hence, we cannot use this vector field to see that $(u,\omega)$ is a strong solution.

To prove $(\z(t), Du(t))=|Du(t)|$ for almost all $t\in (0,T)$, we first fix $t\in(0,T)$ satisfying the previous Steps and observe that we have $\Div\z(t)=\Div\widehat\z(t)$ (by Steps 6, 9 and 13) and $[\z(t),\nu]=[\widehat{\z}(t),\nu]$ (see Step 9). Applying Green's formula,  it yields
\begin{equation*}
  \int_\Omega u(t)\Div\z(t)\, dx+ \int_\Omega (\z(t), Du(t))=  \int_{\partial\Omega} u(t)[\z(t), \nu]\,d\h
\end{equation*}
and
\begin{equation*}
  \int_\Omega u(t)\Div\z(t)\, dx+ \int_\Omega (\widehat\z(t), Du(t))=  \int_{\partial\Omega} u(t)[\z(t), \nu]\,d\h\,,
\end{equation*}
which imply
\begin{equation*}
  \int_\Omega (\z(t), Du(t))=\int_\Omega (\widehat\z(t), Du(t))= \int_\Omega |Du(t)|\,,
\end{equation*}
owed to Steps 10 and 13. Now, it follows from this identity and  $|(\z(t), Du(t))|\le |Du(t)|$ that $(\z(t), Du(t))=|Du(t)|$ as measures. Indeed, take a $|Du|$--measurable set $E\subset\Omega$, then
\begin{multline*}
  \int_\Omega |Du(t)|=\int_\Omega (\z(t), Du(t))=\int_E (\z(t), Du(t))+\int_{\Omega\backslash E} (\z(t), Du(t))\\
  \le\int_E |Du(t)|+\int_{\Omega\backslash E} |Du(t)|=\int_\Omega |Du(t)|
\end{multline*}
and so the inequality becomes equality. Thus, $\int_E (\z(t), Du(t))=\int_E |Du(t)|$.

\medskip
\textbf{STEP 15: Estimates \eqref{est1} and \eqref{est2}.}
\medskip

To check \eqref{est1} we only have to write \eqref{ec-2.1} conveniently as
\begin{align*}
& \bigg(\int_{\partial\Omega}\omega_\e(t,x)^2\,d\h\bigg)^\frac{1}{2}
\\ &\; \le \left(\int_{\partial\Omega}\omega_0(x)^2\,d\h\right)^\frac{1}{2} + \int_{0}^T\left(\int_{\partial\Omega} g_\e(t,x)^2\,d\h\right)^\frac{1}{2}dt \,,
\end{align*}
for all $t\in (0,T)$,
and then apply \eqref{conv-g-L1enL2} and \eqref{conv-omega}.

On the other hand, \eqref{est2} is an easy consequence of \eqref{ec-fron1}.
\end{pf}

\begin{Remark}\rm
It is not difficult to obtain  estimates (other than \eqref{est1} and \eqref{est2}) connecting data and solution, which may have some interest. Indeed, we can easily deduce another estimate starting from \eqref{ec-fron1}:
\begin{align*}
\lambda & \int_\Omega u(t)^2 \,dx + \int_\Omega |Du(t)| +\int_{\partial \Omega} |u(t)| \,d\h
\\ & \le \int_{\partial\Omega}  \big(g(t)-\omega_t(t)\big)\omega(t)\,d\h+\int_{\partial \Omega} |\omega(t)| \,d\h
\\  &  \le\int_{\partial\Omega}  g(t)\omega(t)\,d\h-\frac12\frac d{dt}\int_{\partial\Omega}\omega(t)^2 \,d\h +\int_{\partial \Omega} |\omega(t)| \,d\h   \,.
\end{align*}
Integrating in $[0,t]$ for $t\in (0,T]$, we get
\begin{align*}
\lambda & \int_0^t\int_\Omega u(s)^2 \,dx\, ds + \int_0^t\left[\int_\Omega |Du(s)| +\int_{\partial \Omega} |u(s)| \,d\h\right] \, ds
\\  & +\frac12\int_{\partial\Omega}\omega(t)^2 \,d\h \le\frac12\int_{\partial\Omega}\omega_0^2 \,d\h
\\ & +\int_0^t\int_{\partial\Omega}  g(s)\omega(s)\,d\h\, ds+\int_0^t\int_{\partial \Omega} |\omega(s)| \,d\h \, ds  \,,
\end{align*}
and taking the supremum for $t\in (0,T]$, it yields
\begin{align*}
\lambda & \|u\|_{L^2(0,T;L^2(\Omega))} + \|u\|_{L^1(0,T;BV(\Omega))} +\frac12\|\omega\|_{L^\infty(0,T;L^2(\partial\Omega))}^2
\\ & \le\frac12\|\omega_0\|_{L^2(\partial\Omega)}^2 +\|\omega\|_{L^1(0,T;L^1(\partial\Omega))} +\int_0^T\int_{\partial\Omega}  g(s)\omega(s)\,d\h\, ds
\\ & \le\frac12\|\omega_0\|_{L^2(\partial\Omega)}^2 +\|\omega\|_{L^1(0,T;L^1(\partial\Omega))} +\|\omega\|_{L^\infty(0,T;L^2(\partial\Omega))}\|g\|_{L^1(0,T;L^2(\partial\Omega))}\,.
\end{align*}
Finally, Young's inequality implies
\begin{align*}
\lambda & \|u\|_{L^2(0,T;L^2(\Omega))} + \|u\|_{L^1(0,T;BV(\Omega))}
\\ & \le\frac12\|\omega_0\|_{L^2(\partial\Omega)}^2 +\|\omega\|_{L^1(0,T;L^1(\partial\Omega))}+\frac12\|g\|_{L^1(0,T;L^2(\partial\Omega))}^2\,.
\end{align*}
\end{Remark}

\begin{Remark}\rm
We remark that choosing data in more regular spaces, we get better regularity of the solution. An easy instance is as follows:
If $g \in L^\infty(0, +\infty;L^2(\partial \Omega))$, since the equality $\omega_t(t) =[\z(t),\nu]+g(t)$ holds on $\partial\Omega$, then $\omega_t\in L^\infty(0,+\infty; L^2(\partial\Omega))$ and thus solution $\omega$ is Lipschitz-continuous with respect to the time variable.
\end{Remark}

We finish this section with a Comparison principle and a result on long term behaviour.

\begin{Proposition}
Let $g^1, g^2 \in L^1(0,T;L^2(\partial\Omega))$ and let $\omega_0^1, \omega_0^2 \in L^2(\partial\Omega)$. Denote by $(u^k, \omega^k)$ the strong solution corresponding to data $g^k$ and  $\omega_0^k$, $k=1,2$.

If $g^1(t,x)\le g^2(t,x)$ for almost all $(t, x)\in (0,T)\times \partial\Omega$ and $\omega_0^1(x)\le \omega_0^2(x)$ for almost all $x\in \partial\Omega$, then $\omega^1(t,x)\le \omega^2(t,x)$ for almost all $(t, x)\in (0,T)\times \partial\Omega$ and $u^1(t,x)\le u^2(t,x)$ for almost all $(t,x) \in (0,T)\times \Omega$.
\end{Proposition}
\begin{pf}
It is enough to apply Theorem \ref{compar} having in mind that $\lim\limits_{\e \to 0}u^k_\e(t)=u^k(t)$ in $L^1(\Omega)$ for almost every $t \in (0,T)$.
\end{pf}

\begin{Proposition}
If $g \in L^1(0, +\infty;L^2(\partial \Omega))$ and $\omega_0 \in L^2(\partial \Omega)$, then there exists a sequence $t_n \to +\infty$ and there exist $h\in L^2(\partial \Omega)$ and $v \in L^2(\Omega)\cap BV(\Omega)$ such that $\omega(t_n) \tod h$ weakly in $L^2(\partial\Omega)$, the sequence $\{u(t_n)\}$ converges to $v$ weakly in $L^2(\Omega)$ and strongly in $L^1(\Omega)$ as well as $\{Du(t_n)\}$ converges to $Dv$ $*$-weakly as measures in $\Omega$.
\end{Proposition}
\begin{pf}
Since the datum $g \in L^1(0, +\infty;L^2(\partial \Omega))$ we deduce from estimate \eqref{est1} that
\begin{equation*}
\|\omega\|_{L^\infty (0,+\infty; L^2(\partial\Omega))}\le \|\omega_0\|_{L^2(\partial\Omega)}+\|g\|_{L^1 (0,+\infty; L^2(\partial\Omega))} <+\infty \,.
\end{equation*}
Then, there exist a constant $M>0$ such that $\|\omega(t)\|_{L^2(\partial\Omega)} \le M$ for almost every $t>0$. Therefore, there exist a sequence $t_n \to +\infty$ and a function $h \in L^2(\partial\Omega)$ such that $\omega(t_n) \tod h$ weakly in $L^2(\partial\Omega)$.

On the other hand, from estimate \eqref{est2} we also deduce that
\begin{equation*}
 \lambda\|u(t)\|_{L^2(\Omega)} +\|u(t)\|_{BV(\Omega)} \le M \,.
\end{equation*}
Thus, there exists another sequence $t_n \to +\infty$ and two functions $v_1$ and $v_2$ such that
\begin{equation*}
u(t_n) \tod v_1 \quad \mbox{ in } L^2(\Omega) \,,
\end{equation*}
and
\begin{equation*}
u(t_n) \to v_2 \quad \mbox{ in } L^1(\Omega) \,,
\end{equation*}
with
\begin{equation*}
Du(t_n) \tode Dv_2\quad \mbox{ as measures in } \Omega \,.
\end{equation*}
Finally, due to the uniqueness of the limit, we denote $v=v_1=v_2 \in L^2(\Omega)\cap BV(\Omega)$.
\end{pf}

\section{Continuous dependence on data}

In this section, we get a result which compares solutions of problem \eqref{Dynamic-problem} determined by different data. More precisely, the result allows us to estimate the distance of the solutions depending on the distance of the data.

\begin{Theorem}
Let $(u_1,\omega_1)$ and $(u_2,\omega_2)$ be the strong solution to problem \eqref{Dynamic-problem} with initial data $g_1, g_2 \in L^1(0,T;L^2(\partial\Omega))$ and $\omega_{01}, \omega_{02} \in L^2(\partial \Omega)$ respectively. Then, it holds
\begin{equation*}
 \|\omega_1-\omega_2\|_{L^\infty(0,T;L^2(\partial\Omega))} \le\|\omega_{01}-\omega_{02}\|_{L^2(\partial\Omega)}+\|g_1-g_2\|_{L^1(0,T;L^2(\partial\Omega))}
 \end{equation*}
 and
\begin{equation}\label{uni-2}
\lambda \| u_1-u_2\|_{L^2(0,T;L^2(\Omega))}^2
 \le \dfrac{1}{2} \|\omega_{01}-\omega_{02}\|_{L^2(\partial \Omega)}^2
+ \dfrac{1}{2}\|g_1-g_2\|_{L^1(0,T;L^2(\partial\Omega))}^2 \,.
\end{equation}
\end{Theorem}

\begin{pf}
First, we fix $t \in (0,T)$ such that conditions $(i)$ to $(iv)$ of solution to problem \eqref{Dynamic-problem} hold. Then, we take $u_1(t)-u_2(t)$ as a test function in the condition $(ii)$ corresponding to $(u_1,\omega_1)$. Therefore, using Green's Theorem we get
\begin{align*}
0= & \lambda \int_\Omega u_1(t) (u_1(t)-u_2(t))\,dx + \int_\Omega (\z_1(t),D(u_1(t)-u_2(t)))
\\ & - \int_{\partial \Omega} [\z_1(t),\nu] (u_1(t)-u_2(t)) \,d\h
\,.
\end{align*}

Similarly, we obtain
\begin{align*}
0= &  \lambda\int_\Omega u_2(t) (u_1(t)-u_2(t))\,dx + \int_\Omega(\z_2,D(u_1(t)-u_2(t)))
\\ & - \int_{\partial \Omega} [\z_2(t),\nu] (u_1(t)-u_2(t)) \,d\h \,.
\end{align*}
Now, we combine both equalities to arrive at
\begin{align*}
\lambda & \int_\Omega (u_1(t)-u_2(t))^2 \,dx
\\ &\;+ \int_\Omega \bigg[|Du_1(t)|-(\z_1(t),Du_2(t)) + |Du_2(t)|-(\z_2(t),Du_1(t)) \bigg]
\\ & = \int_{\partial \Omega}  [\z_1(t),\nu] (u_1(t)-u_2(t))  \,d\h + \int_{\partial \Omega}  [\z_2(t),\nu] (u_2(t)-u_1(t)) \,d\h
\\ & =I_1+I_2\,.
\end{align*}
Now, since $(\z_1(t),Du_2(t))\le |Du_2(t)|$ and $(\z_2(t),Du_1(t))\le |Du_1(t)|$, we get the following inequality
\begin{equation}\label{ec-comp}
\lambda\int_\Omega (u_1(t)-u_2(t))^2 \,dx\le I_1+I_2\,.
\end{equation}
We are analyzing $I_1$ and $I_2$. First manipulate $I_1$ using conditions $(iii)$ and $(iv)$:
\begin{align*}
I_1 & =\int_{\partial \Omega} [\z_1(t),\nu] (u_1(t)-\omega_1(t)+\omega_1(t)-\omega_2(t) +\omega_2(t)- u_2(t))\,d\h
\\ & \le - \int_{\partial \Omega} |u_1(t)-\omega_1(t)|\, d\h
\\ & +\int_{\partial \Omega}(g_1(t)-\omega_{1t}(t))(\omega_1(t)-\omega_2(t))\, d\h +\int_{\partial \Omega}|\omega_2(t)- u_2(t)|\,d\h \,.
\end{align*}
In an analogous way we get
\begin{align*}
I_2 & \le - \int_{\partial \Omega} |u_2(t)-\omega_2(t)|\, d\h
\\ & +\int_{\partial \Omega}(g_2(t)-\omega_{2t}(t))(\omega_2(t)-\omega_1(t))\, d\h+\int_{\partial \Omega}|\omega_1(t)- u_1(t)|\,d\h \,,
\end{align*}
and adding both estimates it follows that
\begin{align*}
I_1+I_2 & \le  \int_{\partial \Omega}(g_1(t)-\omega_{1t}(t))(\omega_1(t)-\omega_2(t))\, d\h
\\ & +\int_{\partial \Omega}(g_2(t)-\omega_{2t}(t))(\omega_2(t)-\omega_1(t))\, d\h\,.
\end{align*}
Therefore, \eqref{ec-comp} and H\"older's inequality imply
\begin{align*}
\lambda & \int_\Omega (u_1(t)-u_2(t))^2 \,dx
\\  & \le \int_{\partial \Omega} (g_1(t)-g_2(t)) (\omega_1(t)-\omega_2(t)) \,d\h
\\ & -\int_{\partial \Omega} (\omega_{1t}(t)-\omega_{2t}(t)) (\omega_1(t)-\omega_2(t))  \,d\h
\\  & \le \left(\int_{\partial\Omega} (g_1(t)-g_2(t))^2\,d\h\right)^\frac{1}{2}\left(\int_{\partial\Omega} (\omega_1(t)-\omega_2(t))^2 \,d\h\right)^\frac{1}{2}
\\ & -\int_{\partial \Omega} (\omega_{1t}(t)-\omega_{2t}(t)) (\omega_1(t)-\omega_2(t))  \,d\h\,.
\end{align*}
\\
Moreover, since  $\omega_1, \omega_2\in W^{1,1}(0,T;L^2(\partial \Omega))$,  we know that
\begin{equation*}
\int_{\partial\Omega} (\omega_{1t}(t)-\omega_{2t}(t))(\omega_1(t)-\omega_2(t))\,d\h =\dfrac{1}{2}\dfrac{d}{dt}\int_{\partial\Omega} \big(\omega_{1}(t)-\omega_{2}(t)\big)^2\,d\h \,.
\end{equation*}
Now, let $t \in (0,T)$ and integrate with respect to the time to get
\begin{align*}
\lambda  & \int_0^t \int_\Omega (u_1(s)-u_2(s))^2\,dx\,ds
\\  & \le  \int_0^t \left(\int_{\partial\Omega} (g_1(s)-g_2(s))^2\,d\h\right)^\frac{1}{2}\left(\int_{\partial\Omega} (\omega_1(s)-\omega_2(s))^2 \,d\h\right)^\frac{1}{2} \,ds
\\  & -\dfrac{1}{2} \int_{\partial \Omega}(\omega_{1}(t)-\omega_{2}(t))^2 \,d\h +\dfrac{1}{2} \int_{\partial \Omega}(\omega_{1}(0)-\omega_{2}(0))^2 \,d\h \,.
\end{align*}
So, we have got the main estimate:
\begin{align}\label{uni-3}
2\lambda  & \int_0^t \| u_1(s)-u_2(s)\|_{L^2(\Omega)}^2\,ds+ \|\omega_{1}(t)-\omega_{2}(t)\|_{L^2(\partial\Omega)}^2
\\ \nonumber  & \le 2\int_0^t \|g_1(s)-g_2(s)\|_{L^2(\partial\Omega)}\|\omega_1(s)-\omega_2(s)\|_{L^2(\partial\Omega)} \,ds +\|\omega_{1}(0)-\omega_{2}(0)\|_{L^2(\partial\Omega)}^2 \,,
\end{align}
for all $t\in(0,T)$.
A consequence is the inequality
\begin{align*}
\|\omega_{1}(t)-\omega_{2}(t)\|_{L^2(\partial\Omega)}^2 & \le \|\omega_{1}(0)-\omega_{2}(0)\|_{L^2(\partial\Omega)}^2
\\ & +2\int_0^t \|g_1(s)-g_2(s)\|_{L^2(\partial\Omega)}\|\omega_1(s)-\omega_2(s)\|_{L^2(\partial\Omega)} \,ds \,,
\end{align*}
which, due to an extension of Gronwall's inequality (see \cite{We}), allows us to have
\begin{align*}
\| & \omega_1(t)-\omega_2(t)\|_{L^2(\partial\Omega)} \le \|\omega_{1}(0)-\omega_{2}(0)\|_{L^2(\partial\Omega)}+\int_0^t\|g_1(s)-g_2(s)\|_{L^2(\partial\Omega)}ds
\\  & \le \|\omega_{1}(0)-\omega_{2}(0)\|_{L^2(\partial\Omega)}+\int_0^T\|g_1(s)-g_2(s)\|_{L^2(\partial\Omega)}ds\,,
\end{align*}
for all $t\in (0,T)$. So that
\begin{equation*}
\|\omega_1-\omega_2\|_{L^\infty(0,T;L^2(\partial\Omega))} \le\|\omega_{1}(0)-\omega_{2}(0)\|_{L^2(\partial\Omega)}+\|g_1-g_2\|_{L^1(0,T;L^2(\partial\Omega))}\,.
\end{equation*}
On the other hand, inequality \eqref{uni-3} and Young's inequality imply
\begin{align*}
2\lambda &  \int_0^T \| u_1(s)-u_2(s)\|_{L^2(\Omega)}^2\,ds+\sup_{t\in[0,T]}\|\omega_{1}(t)-\omega_{2}(t)\|_{L^2(\partial\Omega)}^2
\\  & \le ||\omega_{1}(0)-\omega_{2}(0)\|_{L^2(\partial\Omega)}^2
\\ & +2\int_0^T \|g_1(s)-g_2(s)\|_{L^2(\partial\Omega)} \|\omega_{1}(s)-\omega_{2}(s)\|_{L^2(\partial\Omega)} \,ds
\\  & \le ||\omega_{1}(0)-\omega_{2}(0)\|_{L^2(\partial\Omega)}^2
+\|g_1-g_2\|_{L^1(0,T;L^2(\partial\Omega))}^2+ \|\omega_{1}-\omega_{2}\|_{L^\infty(0,T;L^2(\partial\Omega))}^2\,.
\end{align*}
Simplifying, it leads to the desired inequality \eqref{uni-2}.
\end{pf}

\section*{Acknowledgements}
This research has been partially supported  by the Spanish Mi\-nis\-te\-rio de Econom\'{\i}a y Competitividad and FEDER, under project MTM2015--70227--P.
The first author was also supported by Ministerio de Econom\'{\i}a y
Competitividad under grant BES--2013--066655.

\end{document}